\documentclass[letterpaper,journal]{IEEEtran}
\usepackage[]{amsmath,amsfonts}
\usepackage{algorithmic}
\usepackage{array}
\usepackage[caption=false,font=normalsize,labelfont=sf,textfont=sf]{subfig}
\usepackage{textcomp}
\usepackage{stfloats}
\usepackage{url}
\usepackage{verbatim}
\usepackage{graphicx}
\usepackage{cite}
\usepackage{multirow}
\usepackage{siunitx}
\def\BibTeX{{\rm B\kern-.05em{\sc i\kern-.025em b}\kern-.08em
    T\kern-.1667em\lower.7ex\hbox{E}\kern-.125emX}}

\usepackage[inkscapeformat=png]{svg}

\usepackage{mathtools}
\usepackage{lipsum}
\usepackage{amsthm}
\usepackage{flushend}

\theoremstyle{definition}
\newtheorem*{theorem}{Theorem}

\makeatletter
\let\old@ps@headings\ps@headings
\let\old@ps@IEEEtitlepagestyle\ps@IEEEtitlepagestyle
\def\psccfooter#1{%
    \def\ps@headings{%
        \old@ps@headings%
        \def\@oddfoot{\strut\hfill#1\hfill\strut}%
        \def\@evenfoot{\strut\hfill#1\hfill\strut}%
    }%
    \def\ps@IEEEtitlepagestyle{%
        \old@ps@IEEEtitlepagestyle%
        \def\@oddfoot{\strut\hfill#1\hfill\strut}%
        \def\@evenfoot{\strut\hfill#1\hfill\strut}%
    }%
    \ps@headings%
}
\makeatother
\usepackage{balance}
\begin{document}

\title{A Decomposition Method for Solving \\ Sensitivity-Based Distributed Optimal Power Flow}

\author{\IEEEauthorblockN{Mohannad Alkhraijah, 
Devon Sigler, and
Daniel K. Molzahn}
\thanks{\noindent Mohannad Alkhraijah and Devon Sigler are with the Computational Science Center, National Renewable Energy Laboratory, Golden, CO 80401, USA (email: mohannad.alkhraijah@nrel.gov, devon.sigler@nrel.gov). Daniel K. Molzahn is with the School of Electrical and Computer Engineering, Georgia Institute of Technology, Atlanta, GA 30332, USA (email: molzahn@gatech.edu). Support from NSF AI Institute for Advances in Optimization (AI4OPT), \#2112533. Authored in part by the National Renewable Energy Laboratory (NREL), operated by Alliance for Sustainable Energy, LLC, for the U.S. Department of Energy (DOE) under Contract No. DE-AC36-08GO28308. Support from the Laboratory Directed Research and Development (LDRD) Program at NREL. The views expressed in the article do not necessarily represent the views of the DOE or the U.S. Government.}}

\maketitle

\begin{abstract}
Efficiently solving large-scale optimal power flow (OPF) problems is challenging due to the high dimensionality and interconnectivity of modern power systems. Decomposition methods offer a promising solution via partitioning large problems into smaller subproblems that can be solved in parallel, often with local information. These approaches reduce computational burden and improve flexibility by allowing agents to manage their local models. This article introduces a decomposition method that enables a distributed solution to OPF problems. The proposed method solves OPF problems with a sensitivity-based formulation using the alternating direction method of multipliers (ADMM) algorithm. We also propose a distributed method to compute system-wide sensitivities without sharing local parameters. This approach facilitates scalable optimization while satisfying global constraints and limiting data sharing. We demonstrate the effectiveness of the proposed approach using a large set of test systems and compare its performance against existing decomposition methods. The results show that the proposed method significantly outperforms the typical phase-angle formulation with a 14-times faster computation speed on average.
\end{abstract}


\section{Introduction} \label{sec1:introduction}

Electric power systems are undergoing unprecedented changes with a rapid growth of demand and distributed energy resources. System operators face many emerging challenges in coordinating a massive number of generators, maintaining updated systems models, and providing fast solutions to cope with unpredicted changes in generation and demand. This shift motivates research into efficient decentralized operations using distributed optimization and decomposition methods~\cite{molzahn2017survey}. 

Distributed optimization allows multiple agents to collaboratively solve large optimization problems via decomposing the problem into smaller subproblems. Each agent maintains a local model of their subproblem and controls their area of the system. Distributed algorithms can improve scalability and maintainability by allowing local agents to solve local subproblems in parallel. Moreover, decomposing large-scale problems distributes computational burden and memory usage, reducing the amount of stored data at any specific location, and allows local agents to solve local subproblems in parallel. 

\begin{figure}[t]
    \centering
    \includegraphics[width=1\columnwidth]{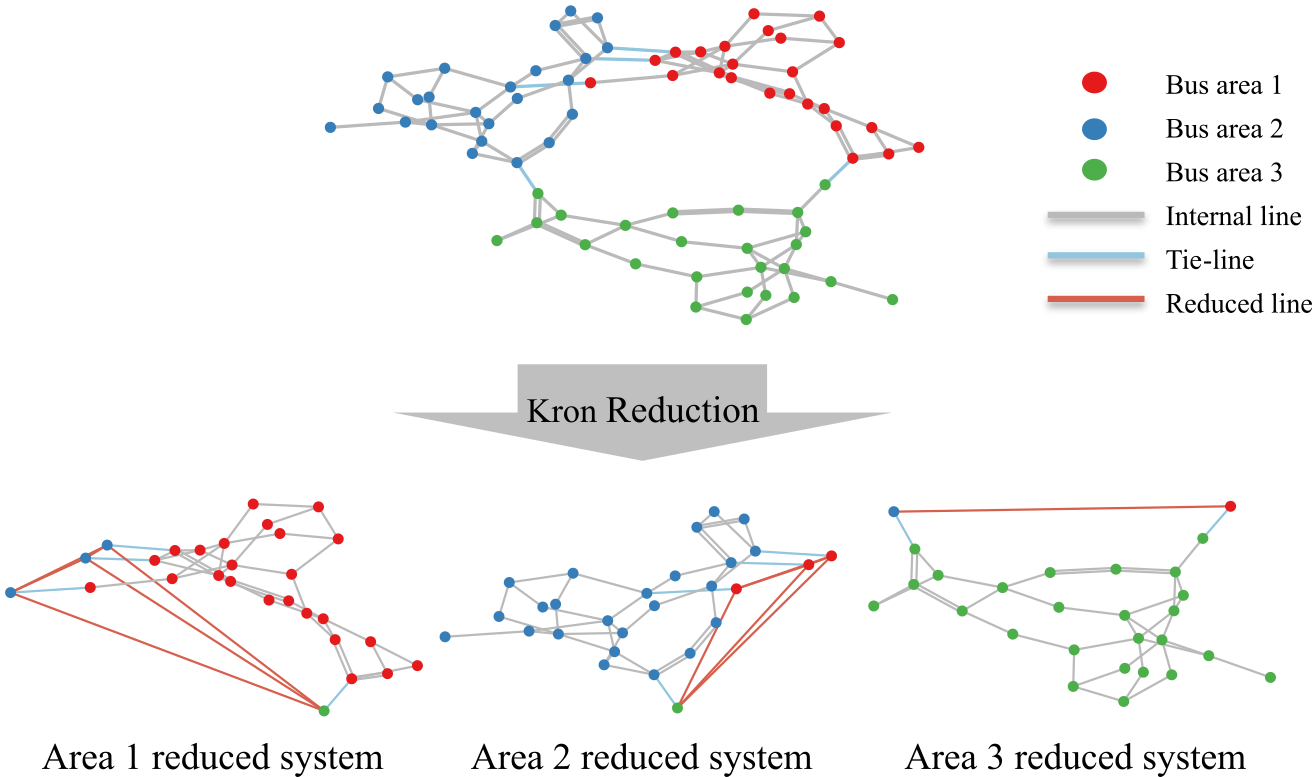}
	\caption[Decomposition diagram]{Decomposition diagram for the IEEE 73-bus Reliability Test System into three areas. The top figure shows the system's three areas with five tie-lines prior the decomposition. The bottom figure shows the reduced system of each area with additional reduced lines between boundary buses.}
	\label{fig:decomposition}
\end{figure}

This article proposes a distributed algorithm that solves optimal power flow (OPF) problems with a novel decomposition method. Unlike typical decomposition methods based on phase-angle formulations, we use Kron reduction to represent the subproblems as shown in Fig.~\ref{fig:decomposition}. The boundary buses of the subproblems capture information from all other subproblems, since each subproblem is a projected version of the original system. We then use the alternating direction method of multipliers (ADMM) algorithm to solve OPF problems with a sensitivity-based formulation. The proposed method converges to the optimal solution and, as the results later in this article show, computationally outperforms typical distributed algorithms based on phase-angle formulations.

\subsection{Related Work}

Several distributed optimization algorithms have been proposed to solve OPF problems, such as ADMM~\cite{7840007}, Analytical Target Cascading~\cite{8994051}, and Auxiliary Problem Principle~\cite{9384964}, and Augmented Lagrangian Alternating Direction Inexact Newton~\cite{8450020}. Distributed algorithms decompose problems into smaller subproblems, typically based on the spatial proximity of the buses. Problem decomposition has been investigated in the literature from two perspectives. The first is how to partition a system, i.e., which variables and constraints belong to which subproblem. The second is how to separate neighboring systems, i.e., what type of consistency constraints should be imposed between coupled subproblems. This article proposes a new approach to separate neighboring systems.

Throughout the literature, we observe three approaches to separate two neighboring areas. The first approach, presented in~\cite{8994051,9384964}, splits the subproblems by duplicating the variables and constraints of the tie-lines connecting two areas. This approach uses auxiliary variables corresponding to the line flow and the voltages at both ends of the tie-lines and imposes consistency constraints. The second approach, used in~\cite{8450020}, decomposes the system by splitting the tie-lines that connect two areas. Then, each area takes a copy of the tie-line with double the admittance value and imposes consistency constraints to equate the voltage values and line flows. The third approach splits the subproblems at the shared buses that belong to two areas~\cite{7840007}. Our previous work in~\cite{Harris2022OnFlow} compares the three approaches and numerically shows that the first method, where the splitting occurs at the tie-lines, achieves the fastest convergence speed compared to the other two approaches. 

Most of the distributed optimization methods proposed in the literature, including all the aforementioned references, use the phase-angle formulation of the OPF problem. An alternative OPF formulation is based on sensitivities, such as Power Transfer Distribution Factors (PTDFs). The PTDF formulation is widely used to solve DC optimal power flow (DCOPF) problems in both research and industry. In particular, PTDFs are often used in market clearing, congestion management, and contingency analysis, as they provide a fast and accurate way to relate power injections and line flows using the DC power flow approximation~\cite{1525117}. The main advantage of using the PTDF formulation is that it reduces the number of decision variables and constraints, which can be significant for large-scale problems~\cite{liu2002effectiveness}. Despite their widespread use in centralized DCOPF formulations, the application of PTDF in distributed OPF settings has received limited attention in the literature.

Reference~\cite{8973905} solves distributed OPF problems with PTDF formulations for market clearing purposes. However, this method uses global PTDFs and the system parameters are known to all participants. Reference~\cite{doi:10.1287/ijoc.2022.0326} proposes a distributed algorithm to solve OPF problems with a local PTDF formulation, and shows that the power flows from neighboring areas are a convex combination of power injections of the external buses. However, this method only allows for two-area decompositions, requiring additional calculations for multi-area decompositions. Another method proposed in~\cite{yang2017decentralized} uses local PTDFs and allows multi-area decompositions. However, this method is only applicable for systems with no loops that pass through multiple areas because there is no guarantee of recovering phase angles from the solution that satisfy the DCOPF constraints. To avoid this,~\cite{zhang2023computational} and~\cite{ding2020zonally} use phase angle constraints for the tie-lines between areas along with local PTDFs. Thus, this approach requires additional equality constraints that limit the benefits of PTDF formulations.

\subsection{Contributions}

This article introduces a method for solving OPF problems with a sensitivity-based formulation using the ADMM algorithm. The main contributions of this article are: 
\begin{itemize}
    \item development of a method based on Kron reduction to decompose OPF problems into multiple subproblems with a mathematical proof guaranteeing consistency with the original central problem,
    \item development of a distributed method to calculate the local PTDF matrices without sharing the local parameters,
    \item proposition of an ADMM-based distributed algorithm for solving OPF problems with the PTDF formulation, and
    \item demonstration of the proposed method using numerical results with 44 test systems that have up to 6500 buses.
\end{itemize}
Unlike the method in~\cite{doi:10.1287/ijoc.2022.0326}, which is only applicable to systems with two areas, the method proposed in this article enables solving OPF problems with any number of areas. Moreover, we propose a systematic method to calculate local PTDF matrices without sharing the local parameters, which involves solving a simple weighted least squares (WLS) problem for each area compared to solving an OPF instance for each column of the local PTDF matrix as proposed in~\cite{doi:10.1287/ijoc.2022.0326}. Finally, the proposed method uses only the PTDF formulation and does not require constraints on phase angles between areas as in~\cite{zhang2023computational} and~\cite{ding2020zonally}. Thus, the proposed method brings the advantages of PTDF-based formulations to distributed settings. To the best of our knowledge, this is the first method that solves OPF problems with the PTDF formulation using distributed optimization that permits multiple-area decompositions. In addition, the reduced representation of the areas is essentially a projection of the full system, which facilitates the exchange of more useful information among areas, leading to a faster convergence speed compared to typical phase-angle formulations, as the results at the end of this article show.

\subsection{Mathematical Notation}

\subsubsection*{Sets and Subsets}
We represent a power system network by a graph $G(\mathcal{N}, \mathcal{E})$, where $\mathcal{N}$ and $\mathcal{E}$ are the sets of buses and lines. We denote the sets of generators and loads as $\mathcal{G}$ and $\mathcal{L}$. We use $\mathcal{G}_n\subset\mathcal{G}$ and $\mathcal{L}_n\subset\mathcal{L}$ to denote the subsets of generators and loads connected to bus $n\in\mathcal{N}$. We use $\mathcal{A}$ to denote the set of areas and define $\mathcal{N}^{a}\subset\mathcal{N}$, $\mathcal{E}^{a}\subset\mathcal{E}$, $\mathcal{G}^{a}\subset\mathcal{G}$, and $\mathcal{L}^{a}\subset\mathcal{L}$ as subsets of buses, lines, generators, and loads belonging to area $a\in\mathcal{A}$. We use $\mathcal{A}_{-a}$ to denote the set of all elements in $\mathcal{A}$ except $a$, i.e., $\mathcal{A}_{-a}\coloneq\mathcal{A}\setminus \{a\}$.

\subsubsection*{Vectors and Matrices}
For a vector $x\in\mathbb{R}^{\mathcal{|\mathcal{N}|}}$ with some set of integers $\mathcal{N}$, the subscript $x_{n}$ denotes the entry corresponds to~$n\in\mathcal{N}$, and the superscript $x^{\alpha}$ denotes the vector containing the entries correspond to the elements in $\alpha\subset\mathcal{N}$. Similarly, for a matrix $A\in\mathbb{R}^{|\mathcal{N}|\times|\mathcal{M}|}$, the subscript $A_{ij}$ denotes the entry corresponds to $i\in\mathcal{N}$ and $j\in\mathcal{M}$, and the superscript $A^{\alpha\beta}$ denotes the block matrix containing the corresponding rows and columns to the elements in $\alpha\subset\mathcal{N}$ and $\beta\subset\mathcal{M}$.

\subsubsection*{Variables and Parameters}
For a power system $G(\mathcal{N}, \mathcal{E})$, the decision variables are the phase angles and power injections of the buses $\theta$ and $p\in\mathbb{R}^{|\mathcal{N}|}$, the power outputs of the generators $g\in\mathbb{R}^{|\mathcal{G}|}$, and the line flows $f\in\mathbb{R}^{|\mathcal{E}|}$. The parameters are the power demands of the loads $d\in\mathbb{R}^{|\mathcal{L}|}$, the upper and lower bounds of the generators $\overline{g}$ and $\underline{g}\in\mathbb{R}^{|\mathcal{G}|}$, the maximum line flows $\overline{f}\in\mathbb{R}^{|\mathcal{E}|}$, and the susceptances of the lines $b\in\mathbb{R}^{|\mathcal{E}|}$. The matrices $B\in \mathbb{R}^{|\mathcal{N}|\times|\mathcal{N}|}$ and $B^E \in \mathbb{R}^{|\mathcal{E}|\times|\mathcal{N}|}$ are the susceptance matrix and the line susceptance matrix, defined in~\cite{molzahn_hiskens-fnt2019}. We use a polynomial cost function $c_k(g_k) = c_{k2} g_{k}^{2} + c_{k1} g_{k} + c_{k0}$ for generator $k\in\mathcal{G}$.

\subsection{Organization}
The remainder of the article is organized as follows. Section~\ref{sec:background} reviews the DCOPF problem with the PTDF formulation. Section~\ref{sec:distributed_opf} presents a distributed optimization algorithm for solving OPF problems with the PTDF formulation. Section~\ref{sec:simulation_results} demonstrates the proposed algorithm via numerical simulations. Section~\ref{sec:conclusions} presents conclusions and future work.

\section{Background} \label{sec:background}

\emph{OPF} is a key optimization problem that finds the operating setpoints while satisfying the power flow equations and engineering constraints. The power flow equations govern the relationship between the voltages, the power injections, and the line flows, while the engineering constraints bound the decision variables. Various OPF formulations proposed in the literature use different linearization and relaxation methods~\cite{molzahn_hiskens-fnt2019}. We consider a linear representation of the OPF problem such as the common DCOPF approximation~\cite{4956966}, which uses the linear DC power flow approximation instead of the nonlinear AC power flow equations. Although we use the DCOPF formulation in this article, the proposed methods apply to other linear OPF representations that allow us to calculate system sensitivities.

\subsection{Power Transfer Distribution Factors}

The \emph{PTDF-OPF} is an equivalent version of the DCOPF problem that uses injection shift factors (ISFs). The \emph{ISFs} calculate the rate of change in line flows with respect to a unit change in power injections balanced by an injection at a reference bus. The \emph{PTDF matrix}, denoted as $H\in\mathbb{R}^{|\mathcal{E}|\times|\mathcal{N}|}$, uses the ISFs to map the power injections to the line flows as
\begin{equation}
    f = H~ p. \nonumber
\end{equation}
Using the DC power flow approximation, we have
\begin{subequations} \label{eq:DC_power_flow} 
\begin{align} 
    & p = B~~\theta, \label{eq:DC_power_balance} \\
    &f = B^E~\theta.
\end{align} \nonumber
\end{subequations} 
\noindent Applying the chain rule on~\eqref{eq:DC_power_flow}, we get the PTDF matrix as 
\begin{equation}\label{eq:PTDF_matrix}
    H = \frac{\partial f}{\partial p} = \frac{\partial f}{\partial \theta}~ \frac{\partial \theta}{\partial p} = B^E~B^{-1}.
\end{equation} 
Since the susceptance matrix $B$ is singular, we invert the reduced matrix defined by omitting the column and row corresponding to a reference bus. This formulation implicitly assumes that the reference bus compensates for any mismatch in power injections. The selection of the reference bus is arbitrary, as we can always compensate for an injection from one bus with an injection of equal magnitude to another bus, which cancels the injection to the reference bus.

\subsection{Central Formulation}

Using the PTDF matrix, the PTDF-OPF formulation is
\begin{subequations} \label{eq:PTDF_opf}
\begin{align}
    & \underset{g, f}{\mbox{minimize}}  && \sum_{k\in \mathcal{G}} c_k(g_{k}), \label{eq:objective} \\
    &\mbox{subject to:} &&  \sum_{k\in \mathcal{G}} g_{k} - \sum_{l\in \mathcal{L}} d_{l} = 0, \label{eq:power_balance} \\
    &&& \!\!\!\!\!\!\!\!\!\! f_e  = \sum_{n\in \mathcal{N}} H_{en} \, (\sum_{k\in \mathcal{G}_n} g_{k}  - \sum_{l\in \mathcal{L}_n} d_{l}), && \forall  e \in \mathcal{E}, \label{eq:line_flow}\\
    &&& \!\!\!\!\!\!\!\!\!\! \underline{g}_{k} \le g_{k} \le \overline{g}_{k}, && \forall  k \in \mathcal{G}, \label{eq:generator_bound} \\
    &&& \!\!\!\!\!\!\!\!\!\! -\overline{f}_e \le f_e \le \overline{f}_e, && \forall  e \in \mathcal{E}. \label{eq:line_flow_bound} 
\end{align}
\end{subequations}
\noindent The objective~\eqref{eq:objective} minimizes the generation cost. Constraints~\eqref{eq:power_balance} and~\eqref{eq:line_flow} are the power balance and power flow equations using the PTDF matrix. Inequalities~\eqref{eq:generator_bound} and~\eqref{eq:line_flow_bound} bound the generators' outputs and the line flows. 

The computational advantages of using PTDF-OPF are significant, especially when a small subset of the line flow bounds are potentially binding. In such cases, we can eliminate the columns and rows of the constraints matrix corresponding to non-binding line flows, which simplifies the problem and enhances computational efficiency. There are similar PTDF-ACOPF formulations that use a first-order Taylor expansion around a nominal operating setpoint to compute linear sensitivities for active and reactive line flows~\cite{wood2013power}.

\section{Distributed Sensitivity-Based \\ Optimal Power Flow} \label{sec:distributed_opf}

This section introduces a distributed PTDF-OPF formulation using the ADMM algorithm. In distributed settings, there are multiple agents, each of which operates a different portion of the power system network, with a fully connected communication network between agents. We first develop a decomposition scheme based on Kron reduction. We then present a distributed method to calculate the local PTDF matrix. After that, we define consistency constraints that ensure recovering a feasible DCOPF solution. Finally, we present the ADMM algorithm to solve distributed PTDF-OPF problems. 

Solving PTDF-OPF problems with distributed optimization requires overcoming two main challenges.
\begin{enumerate}
    \item The value of a line flow depends on the power injections from all buses as implied by~\eqref{eq:line_flow}. To calculate a line flow, we thus need the power injections at the internal buses of other areas, as well as the global PTDF matrix.
    \item The PTDF formulation requires a reference bus that compensates for any mismatch in the system power balance~\eqref{eq:power_balance}. Thus, the power injection of the reference bus depends on the total mismatch in the power balance.
\end{enumerate}
To overcome these challenges, we next present a decomposition scheme that ensures sharing the required information among agents to formulate the distributed PTDF-OPF.

\subsection{System Decomposition} \label{sec:system_decomposition}

Common decompositions to solve OPF problems with distributed algorithms similar to the methods in~\cite{Harris2022OnFlow} are not applicable for PTDF-OPF, as explained previously. Instead, we use Kron reduction to find an equivalent system that eliminates neighboring areas and encompass information from the whole system. Intuitively, the reduced systems are a projected version of the original system with a lower dimension. When multiple agents collaborate in solving the problem, they can recover the central problem with the original dimension.

\emph{Kron reduction} is a well-known technique in electric circuits and power systems analysis to find a reduced equivalent system~\cite{kron1963, dorfler2012kron}. For a power system $G(\mathcal{N}, \mathcal{E})$, let $\alpha \subset \mathcal{N}$ be the subset of buses we want to keep and $\beta \subset \mathcal{N}$ be the subset of the buses that we want to eliminate. We write the power balance equations~\eqref{eq:DC_power_balance} as
\begin{equation}\nonumber
    \begin{bmatrix}
        p^{\alpha} \\ p^{\beta} 
    \end{bmatrix}
    =
    \begin{bmatrix}
        B^{\alpha\alpha} & B^{\alpha\beta} \\ 
        B^{\beta\alpha} & B^{\beta\beta} 
    \end{bmatrix} 
    \begin{bmatrix}
        \theta^{\alpha} \\ \theta^{\beta} 
    \end{bmatrix}.
\end{equation}
\noindent Kron reduction uses Gaussian elimination to get the reduced system. The reduced system's power balance equations are
\begin{equation} \label{eq:kron_reduction}
    p = p^{\alpha} + A~p^{\beta}  = \tilde{B}~\theta^\alpha,
\end{equation}
\noindent where $\tilde{B} = B^{\alpha\alpha} - B^{\alpha\beta}~(B^{\beta\beta})^{-1}~B^{\beta\alpha}$ is the reduced matrix and $A = -B^{\alpha\beta}~(B^{\beta\beta})^{-1}$ is the accompanying matrix.

The distributed PTDF-OPF problem uses the power balance of the reduced system~\eqref{eq:kron_reduction}. The reduced system of any area includes the internal buses, the boundary buses at both ends of tie-lines, and the reference bus. We include the reference bus because the PTDF-OPF uses the reference bus to compensate for changes in power injections.

Comparing the proposed method with the method in~\cite{doi:10.1287/ijoc.2022.0326}, we have three main differences. First, the reference bus is not necessarily internal to each area. This implies that the contribution of the boundary buses to the internal lines is a linear rather than a convex combination of the external buses. Second, the proposed decomposition does not assume that shared buses have zero power injection. Third, instead of solving DCOPF problems to compute the PTDF matrix entries, the proposed decomposition allows for a systematic approach to calculate the local PTDF matrix without sharing the internal parameters via solving a WLS problem, as explained next.

\subsection{Distributed Sensitivity Approximation}

The proposed decomposition requires finding the reduced matrix $\tilde{B}$ and the accompanying matrix $A$ for each area, which involves inverting the matrix $B^{\beta\beta}$ that contains the system parameters from other areas. This section presents a distributed method to calculate the Kron reduction matrices without sharing the internal system parameters.\footnote{The agents do not directly share internal data regarding their areas, which has potential advantages regarding privacy. However, like other distributed optimization methods, there are no rigorous privacy guarantees in the absence of more sophisticated techniques, such as differential privacy methodologies~\cite{9303768}, which can be applied to the method proposed in this article.}

The method proposed in~\cite{doi:10.1287/ijoc.2022.0326} calculates the system sensitivities via solving a sequence of OPF instances. Each OPF instance is used to calculate a single column of a local PTDF matrix. In contrast, our method calculates a local PTDF matrix via solving a single WLS problem that finds the reduced system parameters, which can then be used locally to calculate the PTDF matrix.

Our approach is based on the approximate matrix inversion proposed in~\cite{charalambides2020straggler}. To find the inverse of a square matrix $X\in\mathbb{R}^{n\times n}$ of size $n$, this method solves an unconstrained optimization problem in the form 
\begin{equation} \label{eq:inverse_calculation}
    X^{-1} \approx \underset{Q}{\mbox{arg\,min}}  \;  \left\lVert XQ - I_{n} \right\rVert_{F}^{2},
\end{equation}

\noindent where $\left\lVert \,\cdot\, \right\rVert_{F}$ is the Frobenius norm and $I$ is the identity matrix of appropriate size. Solving~\eqref{eq:inverse_calculation} requires full access to the matrix $X$, which in our case requires full access to the system parameters to solve for the inverse of $B^{\beta\beta}$. We modify this approximation to calculate the Kron reduction for each area without sharing the local system parameters. Instead of inverting $B^{\beta\beta}$, we directly calculate $A = -B^{\alpha\beta}~(B^{\beta\beta})^{-1}$ using distributed optimization.  

Let $\mathcal{A}$ be the set of areas. For area $i\in\mathcal{A}$, we use $\alpha^i$ and $\beta^i$ to denote the set of reduced and eliminated buses, and let $\tilde{B}^i$ and $A^i$ be the reduced and accompanying matrices from the Kron reduction in~\eqref{eq:kron_reduction}. We can approximate $A^i$ as
\begin{equation} \label{eq:PTDF_calculation_central}
    A^{i} \approx \underset{Q}{\mbox{arg\,min}}  \;  \left\lVert Q B^{\beta\beta} + B^{\alpha\beta} \right\rVert_{F}^{2}.
\end{equation}

\noindent We dropped $i$ in $\alpha^{i}$ and $\beta^{i}$ for notational simplicity. 

Let $B^{\beta\beta}_{a}\in\mathbb{R}^{|\beta| \times |\mathcal{N}^a\cap\beta|}$ be a matrix containing the columns of $B^{\beta\beta}$ that correspond to the internal buses of $a\in \mathcal{A}_{-i}$. Similarly, $B^{\alpha\beta}_{a}\in\mathbb{R}^{|\alpha| \times|\mathcal{N}^a\cap\beta|}$ contains the columns of $B^{\alpha\beta}$ that correspond to the internal buses of $a\in \mathcal{A}_{-i}$. We introduce auxiliary variables $Q_{a}$ for all~$a\in \mathcal{A}_{-i}$, and formulate an equivalent problem to~\eqref{eq:PTDF_calculation_central} as 
\begin{equation} \nonumber
     A^{i} \! \approx \! \underset{Q,Q_{a}}{\mbox{arg\,min}} \left \{ \! \sum_{a\in \mathcal{A}_{-i}} \!\!\!   \left\lVert Q_{a} B_{a}^{\beta\beta} \!\! + \! B_{a}^{\alpha\beta} \right\rVert_{F}^{2} \! \mid \! Q = Q_{a}, \forall a \in \mathcal{A}_{-i} \right \}.
\end{equation}
\noindent This is a typical consensus problem that can be solved iteratively using the consensus ADMM algorithm~\cite{boyd2011distributed}. Moreover, we can solve the subproblems analytically as follows
\begin{subequations}
\begin{align}%
    \!\!\!\!\! Q_{a}^{(k+1)} &\!\coloneq \! (B_{a}^{\alpha\beta} (B_{a}^{\beta\beta})^{\mathsf{T}}\!\! -\! \frac{1}{2} \lambda_{a}^{(k)} \!\!  -\! \frac{\rho}{2} Q^{(k)} ) (B_{a}^{\beta\beta}(B_{a}^{\beta\beta})^{\mathsf{T}} \!\! + \! \frac{\rho}{2}I)^{-1}\!, \nonumber \\ &\qquad\qquad\qquad\qquad\qquad\qquad\qquad\qquad \forall a\in\mathcal{A}_{-i}, \\
    \!\!\!\!\! Q^{(k+1)} &\! \coloneq \frac{1}{|\mathcal{A}_{-i}|} \sum_{a\in \mathcal{A}_{-i}} Q_{a}^{(k+1)}, \label{eq:dist_PTDF_para} \\
    \!\!\!\!\! \lambda_{a}^{(k+1)} &\! \coloneq \lambda_{a}^{(k)}\! +\! \rho (Q_{a}^{(k+1)} - Q^{(k+1)}), \qquad\quad \forall a\in\mathcal{A}_{-i},
\end{align}
\end{subequations}
\noindent where $\rho>0$ is a tuning parameter, $\lambda_{a}\in\mathbb{R}^{|\alpha|\times|\beta|}$ is a matrix of Lagrange multipliers, and $I$ is the identity matrix of appropriate size. Agents iteratively solve this problem in parallel, until the values of $\left\lVert Q_{a} - Q \right\rVert^{2}_{F}$ and $\left\lVert Q_{a} B^{\beta\beta}_{a} + B_{a}^{\alpha\beta} \right\rVert^{2}_{F}$, for all $a\in\mathcal{A}_{-i}$, are below a predefined tolerance.

The solution of~\eqref{eq:dist_PTDF_para} converges to the accompanying matrix $A^{i}$, which can be used to calculate the reduced matrix $\tilde{B}^i$ for all $i\in\mathcal{A}$. Finally, agents use the reduced matrix $\tilde{B}^i$ to calculate the local PTDF matrix as described in~\eqref{eq:PTDF_matrix}. 

Using this approach, agents calculate the local PTDF matrices without directly sharing the local parameters. Since the reduced system is a projection of the full system, the local PTDF encapsulates the sensitivities of the full system. In fact, the accompanying matrix $A$ is a linear mapping from the eliminated buses to the boundary buses, which the reduced system uses to calculate the eliminated buses' contributions.

\subsection{Consistency Constraints} \label{sec:consistency_constraints}

This section introduces consistency constraints that ensure recovering the actual problem when using the proposed decomposition. We model the line flow from/to the eliminated buses as fictitious generators (i.e., active power output can take positive and negative values) and impose consistency constraints to equate the reduced generators with the line flows between the areas. Next, we formalize this approach.

\begin{figure}
    \centering
    \includegraphics[width=1\columnwidth]{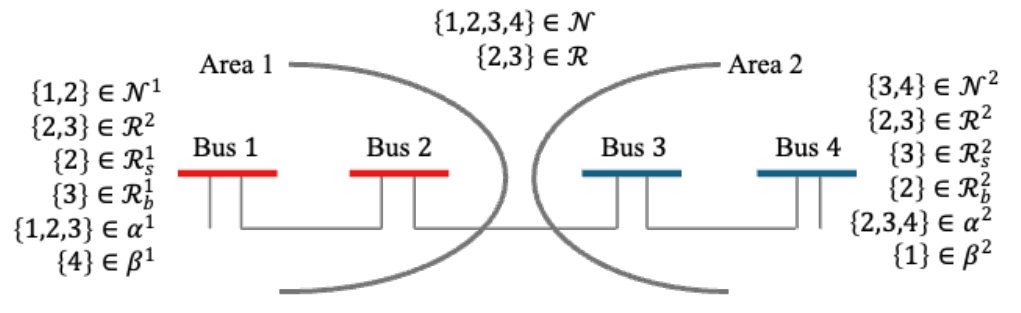}
    \caption{An example of two areas with a tie-line that shows the elements of the sets appearing in the consistency constraints.}
    \label{fig:sets_example}
\end{figure}

Let $\mathcal{R}$ be the set of shared buses, and $\mathcal{R}^i\subset\mathcal{R} $ be the subset of shared buses with area $i$, including the internal and boundary buses denoted $\mathcal{R}^{i}_{s}\subset \mathcal{R}^i$ and $\mathcal{R}^{i}_{b}\subset \mathcal{R}^i$, respectively. The subproblem of area $i$ includes reduced generators connected to boundary buses $\mathcal{R}^i_{b}$ that represent the line flows from neighboring agents. Recall that $\mathcal{N}^i$ is the set of internal buses in area $i$. Thus, $\alpha^{i} \coloneq \mathcal{N}^{i} \cup \mathcal{R}^{i}_{b}$ and $\beta^{i} \coloneq \mathcal{N} \setminus \alpha^{i}$. An example of the sets defined in this section is shown in Fig.~\ref{fig:sets_example}. Using this notation, we define the consistency constraints as
\begin{multline} \label{eq:dist_consistency}
p^{i}_{n} + \sum_{v\in\beta^{j}\cap\mathcal{N}^{i}} \!\!\!\! A^{j}_{nv}~p^{i}_{v} 
- p^{j}_{n} - \sum_{v\in\beta^{i}\cap\mathcal{N}^{j}} \!\!\!\! A^{i}_{nv}~p^{j}_{v} \\
- \sum_{a\in\mathcal{A}_{-ij}} \sum_{v\in\beta^{i}\cap\mathcal{N}^{a}} \!\!\!\! A^{i}_{nv}~p^{a}_{v} 
+ \sum_{a\in\mathcal{A}_{-ij}} \sum_{v\in\beta^{j}\cap\mathcal{N}^{a}} \!\!\!\! A^{j}_{nv}~p^{a}_{v}=0.
\end{multline}
Imposing the consistency constraints~\eqref{eq:dist_consistency} for all shared variables $n\in\mathcal{R}$ ensures the decomposed problem is consistent with the central PTDF-OPF~\eqref{eq:PTDF_opf}, as stated in the theorem below. 
\begin{theorem} \label{th:recovering_phase_angles}
Given a connected system $G(\mathcal{N},\mathcal{E})$ with a susceptance matrix $B$. Assume the susceptances are independent and have non-zero values. Using the decomposition described in Section~\ref{sec:system_decomposition} to decompose the problem into subproblems $\mathcal{A}$, a solution $p^i$ and $\theta^i\in\mathbb{R}^{|\mathcal{N}^i|}$, $i\in\mathcal{A}$, is a solution to the power balance equations of the reduced systems~\eqref{eq:kron_reduction} and satisfies the consistency constraints~\eqref{eq:dist_consistency} if and only if there is a unique solution $p$ and $\theta\in\mathbb{R}^{|\mathcal{N}|}$ to the power balance equations of the original system~\eqref{eq:DC_power_flow} such that $\theta_n = \theta_n^i$ for all $n\in\mathcal{N}^i, i\in\mathcal{A}$.
\end{theorem}
The proof is shown in Appendix~\ref{ap:proof}

\subsection{Distributed Formulation}

This section describes the use of the ADMM algorithm to solve PTDF-OPF with the decomposition described in Section~\ref{sec:system_decomposition}. We introduce auxiliary variables $z^a \in \mathbb{R}^{|\mathcal{R}|}$, for all $a\in\mathcal{A}$, and a consistency constraint~\eqref{eq:dist_consistency} for all shared buses $n\in \mathcal{R}$. We rewrite the consistency constraint of a shared bus $n$ inside area $i$ and on the boundary of area $j$ as
\begin{subequations} \label{eq:consistency_admm}
\begin{align}
  & z^{i}_{n} = p^{i}_{n} + \sum_{v\in\beta^{j}\cap\mathcal{N}^{i}} A^{j}_{nv}~p^{i}_{v}, \label{eq:inside_bus}\\
  & z^{j}_{n} = - p^{j}_{n} - \sum_{v\in\beta^{i}\cap\mathcal{N}^{j}} A^{i}_{nv}~p^{j}_{v}, \label{eq:boundary_bus} \\
  & z^{a}_{n} = - \!\!\!\!\! \sum_{v\in\beta^{i}\cap\mathcal{N}^{a}} \!\!\! A^{i}_{nv}~p^{a}_{v} +  \sum_{v\in\beta^{j}\cap\mathcal{N}^{a}} \!\!\! A^{j}_{nv}~p^{a}_{v}, \forall a\in \mathcal{A} _{-ij}, \label{eq:other_shared_bus}\\
  & \sum_{a\in\mathcal{A}} z^{a}_{n} = 0. \label{eq:central_consistency_constraints}
\end{align}
\end{subequations}
Each agent uses one of~\eqref{eq:inside_bus}--\eqref{eq:other_shared_bus} in their subproblem, while~\eqref{eq:central_consistency_constraints} is a global constraint that ensures satisfying~\eqref{eq:dist_consistency}. For the example above with the shared bus $n$ inside area $i$ and on the boundary of area $j$, agent $i$ takes~\eqref{eq:inside_bus}, agent $j$ takes~\eqref{eq:boundary_bus}, and the other agents take~\eqref{eq:other_shared_bus}. We collect the right-hand side of~\eqref{eq:consistency_admm} in a vector $\zeta^i\in\mathbb{R}^{|\mathcal{R}|}$ that contains the shared variables associated with agent~$i$, including the shared variables that are not in the reduced system of agent $i$, such that
\begin{equation} \label{eq:consistency_area_admm}
    \!\!\! \zeta^{i}_{n} \! \coloneqq \!\!
    \begin{cases}
       p^{i}_{n} + \underset{v\in\beta^{j}\cap\mathcal{N}^{i}}{\sum}  A^{j}_{nv}~p^{i}_{v}, 
       &~\!\!\!\!\mbox{if}~n\in \mathcal{R}^{i}_{s}~\cap~\mathcal{R}^{j}_{b}, \\
       - p^{i}_{n} - \underset{v\in\beta^{j}\cap\mathcal{N}^{i}}{\sum} A^{j}_{nv}~ p^{i}_{v}, 
       &~\!\!\!\!\mbox{if}~n\in \mathcal{R}^{j}_{s}~\cap~\mathcal{R}^{i}_{b}, \\
       - \!\!\!\!\! \underset{v\in\beta^{a}\cap\mathcal{N}^{i}}{\sum} \!\!\!\!\!\! A^{a}_{nv}~p^{i}_{v} + \!\!\!\! \underset{v\in\beta^{j}\cap\mathcal{N}^{i}}{\sum} \!\!\!\!\!\! A^{j}_{nv}~p^{i}_{v}, 
       &~\!\!\!\!\mbox{if}~n\in \mathcal{R}^{a}_{s}~\cap~\mathcal{R}^{j}_{b}.
    \end{cases}
\end{equation}
\noindent Each agent uses one of the three expressions in~\eqref{eq:consistency_area_admm} based on where bus~$n$ is located and shared with whom. Note that the expression of~$\zeta^i$ includes power injections in area $i$ defined as
\begin{equation} \label{eq:consistency_area_admm_variables}
    p^{i}_{n} \coloneqq \sum_{k\in \mathcal{G}^{i}_{n}} g_{k} - \sum_{l\in \mathcal{L}^{i}_{n}} d_{l}, \; \forall n \in \mathcal{N}^i\cup\mathcal{R}_{b}^{i}.
\end{equation}
\noindent Recall that $\mathcal{G}^i_{n}$ and $\mathcal{L}^{i}_{n}$ denote the subsets of generators and loads connected to bus~$n$ in area~$i$, and $\mathcal{R}_{b}^{i}$ is the set of boundary buses of area~$i$, which have reduced generators.

We then iteratively solve the PTDF-OPF problem using the ADMM algorithm as follows
\begin{subequations}
\begin{align}
    & (p^{i}, g^{i}, f^{i},\zeta^{i})^{(k+1)}\!\! \coloneq \! \underset{p, g, f,\zeta}{\mbox{arg\,min}} \!  \left \{ \sum_{k\in \mathcal{G}^{i}} \!\! c_{k} (g_{k}) + \frac{\rho}{2}  \left\lVert z^{i(k)} \! - \zeta \right\rVert_{2}^{2} \!\!\! \right. \nonumber \\ & \left. \!\!\!  + (y^{i(k)})^{\mathsf{T}} \! \left (z^{i(k)} \! - \zeta \right ) \! \mid \eqref{eq:power_balance} \mbox{--} \eqref{eq:line_flow_bound}, \eqref{eq:consistency_area_admm}, \eqref{eq:consistency_area_admm_variables} \right \}, \;\, \forall i \in \mathcal{A}, \label{eq:dist_PTDF_opf} \\
   & z^{i(k+1)} \coloneq \zeta^{i(k+1)} - \frac{1}{|\mathcal{A}|} \sum_{j\in \mathcal{A}} \zeta^{j(k+1)}, \qquad\quad\;\forall i \in \mathcal{A}, \label{eq:central_dist_PTDF_opf} \\
   &y^{i(k+1)} \coloneq y^{i(k)} + \rho \left ( z^{i(k+1)} - \zeta^{i(k+1)} \right ),\quad\;\; \forall i \in \mathcal{A}, \label{eq:dual_dist_PTDF_opf}
\end{align}
\end{subequations}
\noindent where $\rho>0$ is a tuning parameter and $y^{i(k)}\in\mathbb{R}^{|\mathcal{R}|}$ is a vector of Lagrange multipliers. The value of $z^{i(k)}$ is a projection of the shared variables $\zeta^{i(k)}$ into the linear subspace given by~\eqref{eq:central_consistency_constraints}, which we analytically compute using the derivation in~\cite[sec~5]{eckstein2015understanding}.
\begin{figure}[t]
    \centering
    \includegraphics[width=1\columnwidth]{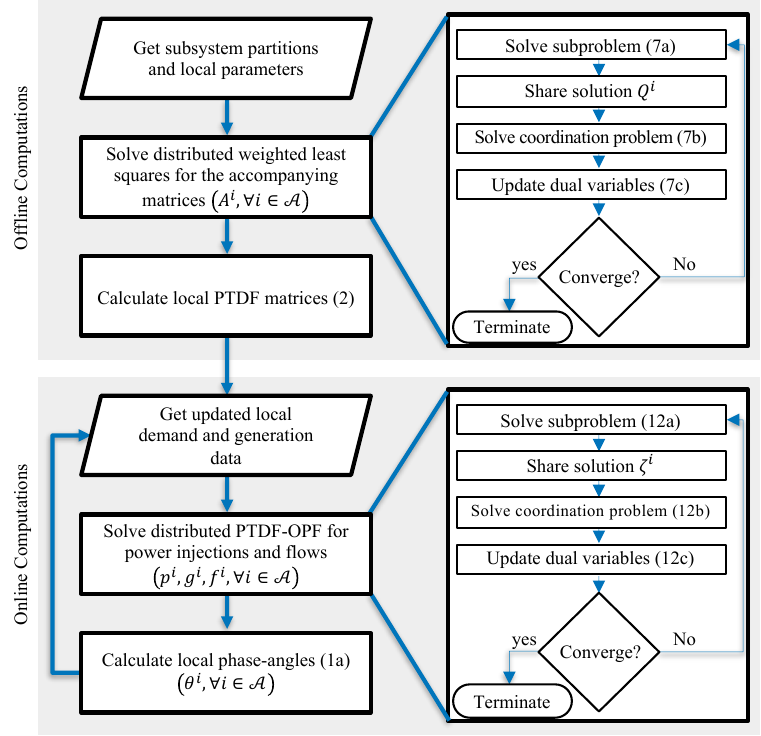}
    \caption{Flowchart of the proposed method. In the top gray block, agents find the reduced systems parameters offline. In the bottom block, agents solve distributed PTDF-OPF to find the operational setpoints online.}
    \label{fig:flowchart}
\end{figure}
The agents then repeat solving the subproblem~\eqref{eq:dist_PTDF_opf}, share the solutions of the auxiliary variables $\zeta^i$, solve the coordination problem~\eqref{eq:central_dist_PTDF_opf}, and update the Lagrange multipliers $y^i$ using~\eqref{eq:dual_dist_PTDF_opf}, until the solution converges. The algorithm converges when the $l_2$-norm of the primal residual $\left\lVert z^{i(k)} - \zeta^{i(k)} \right\rVert_{2}^{2}$ and the dual residual $\left\lVert \rho (z^{i(k)} - z^{i(k-1)}) \right\rVert_{2}^{2}$ for all $i\in\mathcal{A}$ are below a predefined tolerance. Since the ADMM algorithm converges to the optimal solution for convex problems~\cite{boyd2011distributed} and the distributed PTDF-OPF subproblems~\eqref{eq:dist_PTDF_opf}--\eqref{eq:dual_dist_PTDF_opf} are convex, the algorithm converges to the optimal solution.

The shared data at the end of each iteration in typical distributed phase-angle formulations consists of local information about the boundary buses. In contrast, the data shared at the end of each iteration of the proposed method encapsulate information from the entire system. This has a significant impact on the convergence speed of distributed optimization, as the results in the next section show. Moreover, any feasible solution of a subproblem can be used to obtain a solution to the power balance equations of the full system~\eqref{eq:DC_power_flow}.

The overall flowchart of the proposed method is shown in Fig.~\ref{fig:flowchart}. In summary, the proposed method solves PTDF-OPF problems~\eqref{eq:PTDF_opf} using the ADMM algorithm. The method employs Kron reduction to decompose the systems into multiple areas and solves WLS problems to calculate the local PTDF matrices without sharing the local parameters. The local PTDF calculation step is an offline process that agents need to perform once or when the system parameters change. Finally, agents use the local PTDF matrix to solve distributed PTDF-OPF. To the best of our knowledge, this is the first method that solves PTDF-OPF problems using distributed optimization with multiple-area decomposition.

\section{Simulation Results} \label{sec:simulation_results}

This section presents simulation results for solving OPF problems with the PTDF-OPF formulation using the ADMM algorithm and a comparison with the phase-angle formulation. We consider 44 test systems from the PGLib-OPF library~\cite{PGLib} with sizes ranging from 50 to 6500 buses. We decompose the test systems based on the area of the buses. For test systems with a single area, we partition the systems into five areas using the KaHyPar partition algorithm~\cite{10.1145/3529090}. 
\begin{figure}[t]
    \centering
        \includegraphics[width=1\columnwidth]{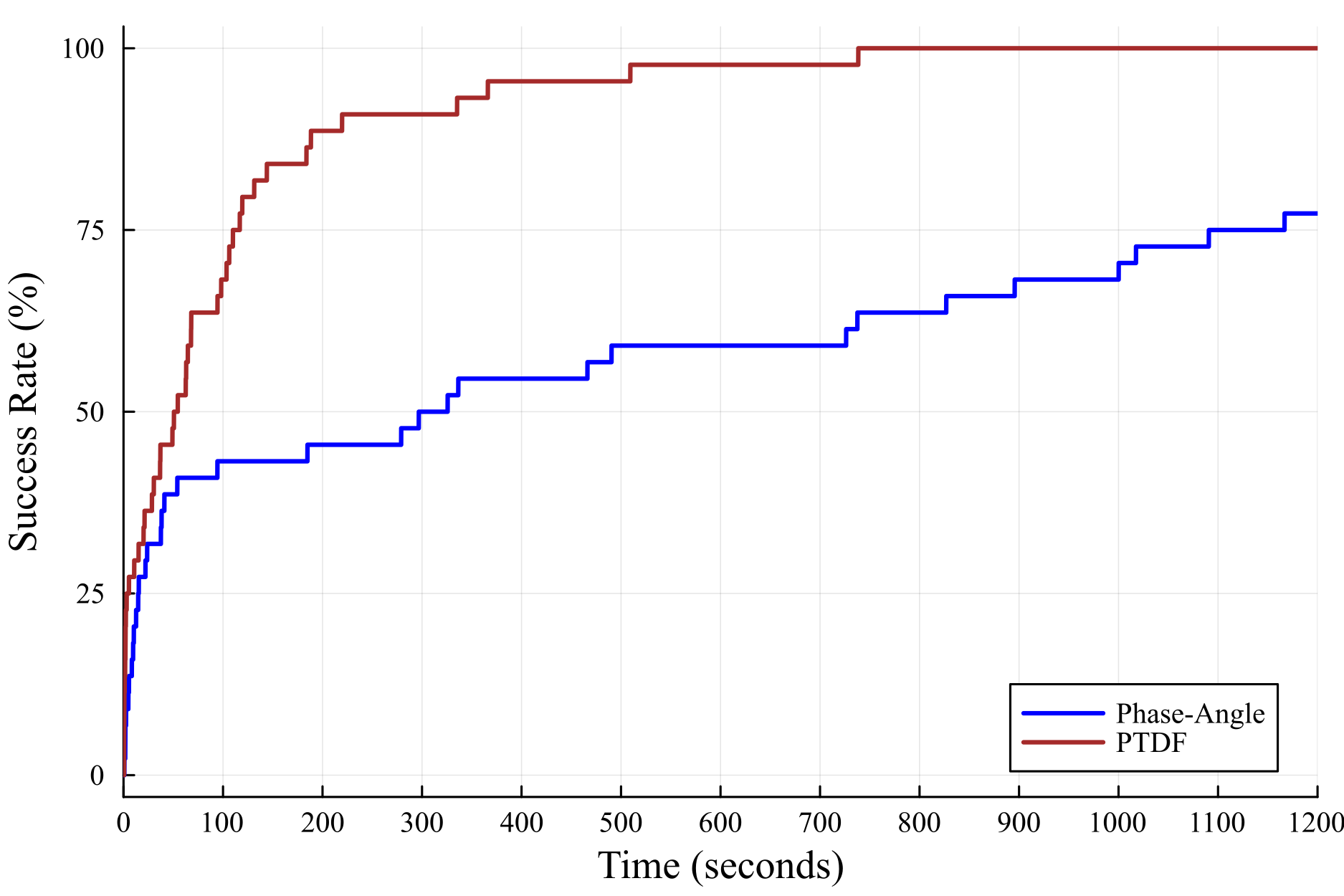}
	\caption[successtime]{Success rate over convergence time for the phase-angle and the proposed PTDF formulations using 44 test cases from PGlib.}
	\label{fig:succes_vs_time}
\end{figure}
The results use PowerModels~\cite{8442948} and \mbox{PowerModelsADA} \cite{alkhraijah2023powermodelsada} libraries in the Julia programming language to solve central and distributed OPF problems and use the Ipopt solver~\cite{Andreas_ipopt}. We produced the results using a computer with a 4.5~GHz, 14-core CPU and 48~GB of memory. The simulation results in this section use a parallel implementation of the ADMM algorithm.

We set the parameter of the ADMM algorithm $\rho = 1000$ and choose primal and dual termination tolerances equal to $1\times10^{-3}$. We consider the algorithm to be successful in solving an OPF instance if it converges to the optimal solution in less than 3600 seconds. We use a lazy constraint method with the proposed PTDF method, where we add line flow variables and constraints if the corresponding line flow bounds are violated in the previous iteration.  

The success rate over time of the PTDF and phase-angle formulations is shown in Fig.~\ref{fig:succes_vs_time}. The proposed method with the PTDF formulation successfully converged to the optimal solution in all 44 test cases in less than 800 seconds, while the phase-angle formulation failed in two test cases. Moreover, the proposed PTDF formulation converged in 75\% of the test cases in around 100 seconds, while the phase-angle formulation took more than 1000 seconds to converge in 75\% of the test cases. When comparing the number of iterations, the difference is very evident. The PTDF formulation required 1320 iterations to converge in all test cases, while the phase angle formulation required over 10,000 iterations to converge in 50\% of the test cases.

\begin{figure}[t]
    \centering
        \includegraphics[width=1\columnwidth]{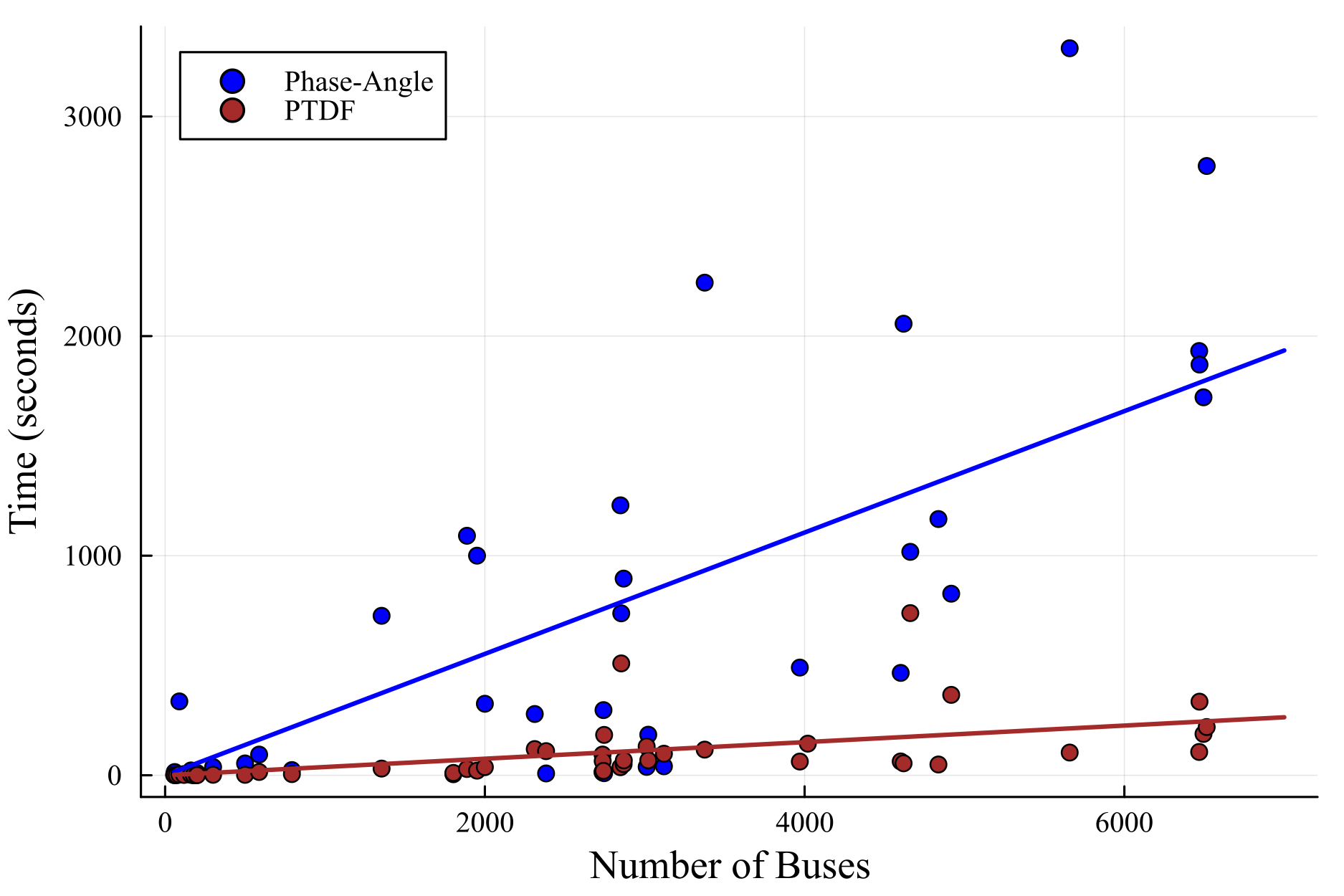}
	\caption[timebuses]{Convergence time with respect to number of buses for 44 test cases. The dots represent instances, and the lines represent approximate trends of the convergence time as the number of iterations increases.}
	\label{fig:time_vs_buses}
\end{figure}
We also observe that, on average, the computation time of a single iteration in the proposed PTDF method is four times greater than the phase-angle formulation, mainly because the PTDF method checks for constraint violations at each iteration. However, the average number of iterations in the proposed PTDF is 50 times less, significantly reducing the overall computation time.

\begin{figure}[t]
    \centering
        \includegraphics[width=1\columnwidth]{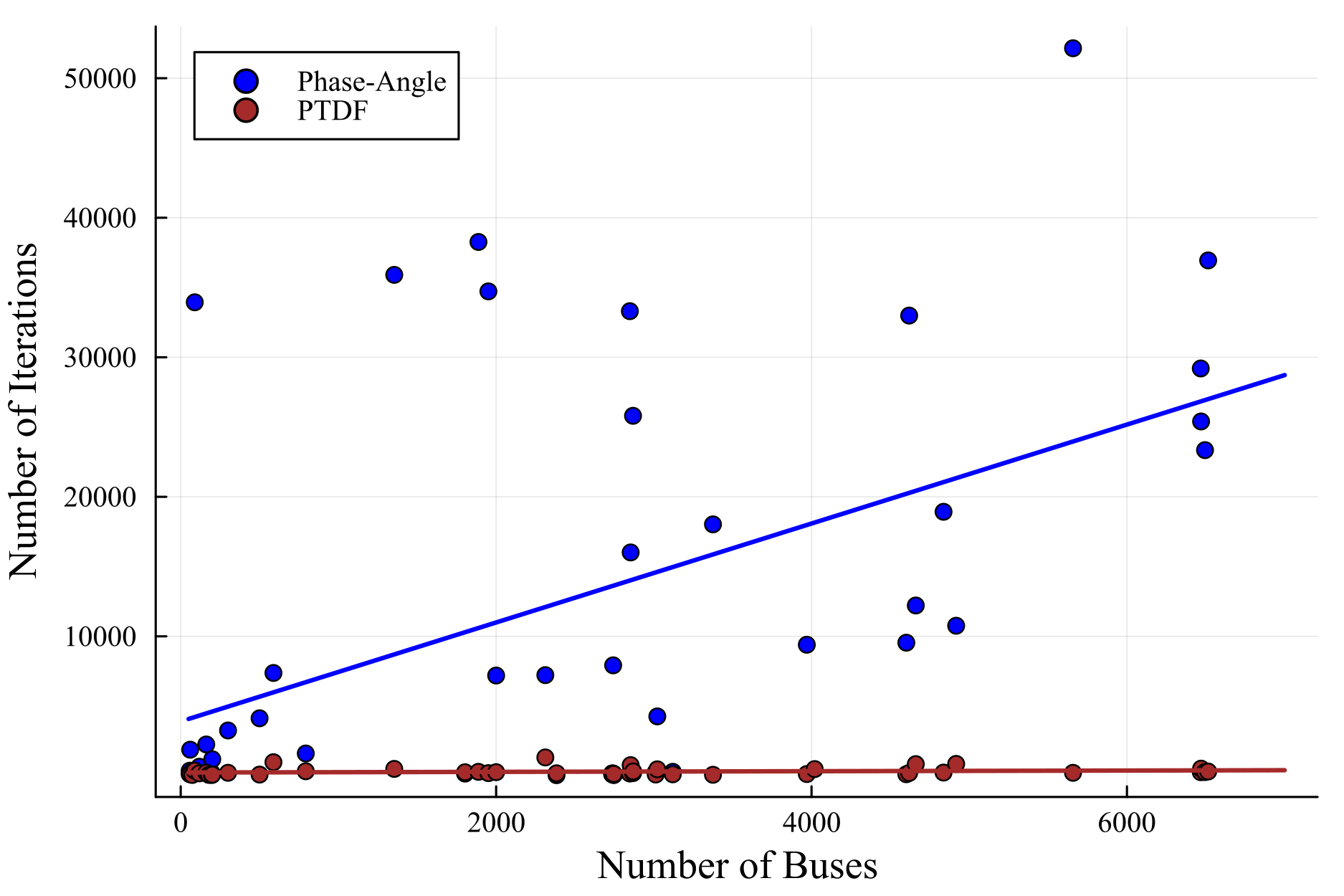}
	\caption[iterationbuses]{Convergence iteration with respect to number of buses for 44 test cases. The dots represent instances, and the lines are approximate trends of the number of iterations required to converge as the number of iterations~increases.}
    \label{fig:iteration_vs_buses}
\end{figure}

We also investigated the convergence time and iterations with respect to the size of the system, as shown in Fig.~\ref{fig:time_vs_buses} and Fig.~\ref{fig:iteration_vs_buses}. The figures show when the algorithms successfully solved the 44 test cases and the approximated trend as the number of buses increases. We observe that as the size of the system increases, both formulations require a longer time to converge. However, the average number of iterations required for the PTDF formulation to converge is almost constant as the size of the system increases, as shown by the red line in Fig.~\ref{fig:iteration_vs_buses}. The list of test systems and the results for the PTDF and phase-angle formulations are given in Table~\ref{tab:PTDF_comparison} in Appendix~\ref{ap:results}.

\section{Conclusion} \label{sec:conclusions}
Decomposition methods provide many advantages in solving large-scale OPF problems, as they distribute computational burden and storage usage among multiple agents, in addition to allowing parallel computation and limiting data sharing. However, typical distributed formulations of the OPF problem suffer from a very slow convergence rate, hindering their usage in practice. This article proposes a new decomposition method based on the Kron reduction to solve distributed OPF problems with the PTDF formulation. The numerical results show that the proposed method significantly outperforms the phase-angle formulation in terms of convergence speed and consistency. 

Sensitivity-based formulations, such as the PTDF-OPF, have been widely used in practice for contingency analysis and congestion management, particularly in deregulated markets. Traditionally, system operators conduct contingency analysis and congestion management independently of neighboring systems. Using PTDF formulations with distributed algorithms to analyze interconnected systems offers a promising application of the proposed method. Furthermore, these methods can facilitate the clearing of interconnected markets and the valuation of energy exchange costs with minimal data sharing.

\appendices

\section{Proof of Theorem~\ref{th:recovering_phase_angles}} \label{ap:proof}

\begin{proof}

We prove the forward implication by induction on the number of areas. The implication is trivially true for $|\mathcal{A}|=1$. Assuming it is true for $|\mathcal{A}|=k$, we need to show that it is also true for $|\mathcal{A}|=k+1$.
It is sufficient to show that the implication is true for two areas, since we can divide the areas into two groups $\mathcal{A}_1$ and $\mathcal{A}_2$, with $|\mathcal{A}_1|=k$ and $|\mathcal{A}_2|=1$, and the induction hypothesis implies the implication for $\mathcal{A}_1$.

Consider a system of two areas $\mathcal{A} =\{1,2 \}$. Recall that $\mathcal{N}^i$ is the set of internal buses and $\mathcal{R}^i$ is the set of shared buses, which consists of the sets of internal and boundary shared buses $\mathcal{R}^i_s$ and $\mathcal{R}^i_b$. Let $\mathcal{N}^{i}_{t}$ be the set of internal buses that are not shared, i.e., $\mathcal{N}^{i}_{t} := \mathcal{N}^{i} \setminus\mathcal{R}^{i}_s$.  We use $p_1$, $p_s$, $p_b$, and $p_2$ to denote the power injections of the buses in $\mathcal{N}^{1}_{t}$, $\mathcal{R}^{1}_{s}$, $\mathcal{R}^{1}_{b}$, and $\mathcal{N}^{2}_{t}$, respectively. Similarly, we use $p^{i}_{t}$, $p^{i}_{s}$, and $p^{i}_{b}$ to denote the vectors of area $i$ reduced power injections at the buses in $\mathcal{N}^{i}_{t}$, $\mathcal{R}^{i}_{s}$, and $\mathcal{R}^{i}_{b}$, respectively. Note that $\mathcal{R}^{1}_{s} = \mathcal{R}^{2}_{b}$ and  $\mathcal{R}^{1}_{b} = \mathcal{R}^{2}_{s}$. The power balance equations prior to Kron reduction~\eqref{eq:DC_power_balance} are
\begin{equation}\label{eq:global_power_flow_1}
    \begin{bmatrix}
        p_{1} \\
        p_{s} \\
        p_{b} \\
        p_{2}
    \end{bmatrix}
    = 
    \begin{bmatrix}
        B_{11}  & B_{1s}    & 0        & 0\\
        B_{s1}  & B_{ss}    & B_{sb}   & 0\\
        0       & B_{bs}    & B_{bb}   & B_{b2}\\
        0       & 0         & B_{2b}   & B_{22}\\
    \end{bmatrix} 
    \begin{bmatrix}
        \theta_{1} \\
        \theta_{s} \\
        \theta_{b} \\
        \theta_{2}
    \end{bmatrix}.
\end{equation}
\noindent Thus, the reduced power balance of area 1 and 2 are
\begin{align}
    \begin{bmatrix}
        p^{1}_{t} \\
        p^{1}_{s} \\
        p^{1}_{b}
    \end{bmatrix}
     \!\! = \!\! 
    \begin{bmatrix}
        p_{1} \\
        p_{s} \\
        p_{b} -  B_{b2} B_{22}^{-1} p_{2}\\
    \end{bmatrix}
    \!\! = \!\! 
    \begin{bmatrix} \label{eq:area1_power_flow}
        B_{11}  & B_{1s}    & 0 \\
        B_{s1}  & B_{ss}    & B_{sb} \\
        0       & B_{bs}    & B_{22}^{'} \\
    \end{bmatrix} \!\!\!
    \begin{bmatrix}
        \theta^{1}_{t} \\
        \theta^{1}_{s} \\
        \theta^{1}_{b} \\
    \end{bmatrix},
\end{align}
\begin{align}
    \begin{bmatrix}
        p^{2}_{b} \\
        p^{2}_{s} \\
        p^{2}_{t}
    \end{bmatrix}
    \!\! = \!\! 
    \begin{bmatrix}
        p_{s} -  B_{s1} B_{11}^{-1} p_{1}\\
        p_{b} \\
        p_{2} \\
    \end{bmatrix}
    \!\! = \!\! 
    \begin{bmatrix}\label{eq:area2_power_flow}
        B_{11}^{'}  & B_{sb}    & 0 \\
        B_{bs}      & B_{bb}    & B_{b2}\\
        0           & B_{2b}    & B_{22}\\
    \end{bmatrix} \!\!\!
    \begin{bmatrix}
        \theta^{2}_{b} \\
        \theta^{2}_{s} \\
        \theta^{2}_{t} \\
    \end{bmatrix},
\end{align}
\noindent where $B_{11}^{'} = B_{ss}\! - \! B_{s1} B_{11}^{-1} B_{1s}$ and $B_{22}^{'} = B_{bb}\! - \! B_{b2} B_{22}^{-1} B_{2b}$. Finally, we write the consistency constraints~\eqref{eq:dist_consistency} as
\begin{align}
p^{1}_{b} &= p_{s}^{2} - B_{b2} B_{22}^{-1} p^{2}_{t}, \label{eq:proof_consistency_1} \\
p^{2}_{b} &= p_{s}^{1} - B_{s1} B_{11}^{-1} p^{1}_{t}. \label{eq:proof_consistency_2}
\end{align}

We want to show that, given a solution of power injections from the reduced subproblems that satisfy the consistency constraints, we can recover a unique solution to the original system such that the phase angles of the buses that appear in the two reduced systems are the same. Formally, if $p^1_t, p^1_s, p^1_b$ and $p^2_t, p^2_s, p^2_b$ satisfy the reduced power balance~\eqref{eq:area1_power_flow} and~\eqref{eq:area2_power_flow}, and the consistency constraints~\eqref{eq:proof_consistency_1} and~\eqref{eq:proof_consistency_2}, then the phase angles of the shared buses $\theta^{1}_{s}=\theta^{2}_{b}$ and $\theta^{1}_{b}=\theta^{2}_{s}$.

Using the consistency constraints~\eqref{eq:proof_consistency_1}, we have 
\begin{align}
    p^{1}_{b}   & = p_{s}^{2} - B_{b2} B_{22}^{-1} p^{2}_{t} \nonumber\\ 
                & =  B_{bs} \theta^{2}_{b} + B_{bb} \theta^{2}_{s} +  B_{b2} \theta^{2}_{t} - B_{b2} B_{22}^{-1} (B_{2b}\theta^{2}_{s} + B_{22} \theta^{2}_{t})  \nonumber\\
                & = B_{bs} \theta^{2}_{b} + B_{22}^{'} \theta^{2}_{s}.\label{eq:area1_boundary}
\end{align}
\noindent We use~\eqref{eq:area2_power_flow} in the first equality and rearrange the terms in the second. We also have from~\eqref{eq:area2_power_flow} that
\begin{align}
    p^{2}_{b} & = B_{11}^{'} \theta^{2}_{b} + B_{sb} \theta^{2}_{s} .\label{eq:area2_boundary_term}
\end{align}
The system of~\eqref{eq:area1_boundary} and~\eqref{eq:area2_boundary_term} consists of $|\mathcal{R}|$ unknowns and $|\mathcal{R}|$ linearly independent equations. The linear independency follows from the assumption that line susceptances are independent~\cite[Th.~1]{9889186}. Thus, the system has a unique solution. Similarly, we can show that
\begin{align}
    p^{1}_{b} & = B_{bs} \theta^{1}_{s} + B_{22}^{'} \theta^{1}_{b}, \nonumber\\
    p^{2}_{b} & = B_{11}^{'} \theta^{1}_{s} + B_{sb} \theta^{1}_{b}. \nonumber
\end{align}
Since both systems are the same and have a unique solution, their solution must be the same, i.e., $\theta^{1}_{s} = \theta^{2}_{b}$ and $\theta^{1}_{b} = \theta^{2}_{s}$. Thus, we can recover a unique solution to~\eqref{eq:global_power_flow_1} by setting the value of $\theta_{1} = \theta^{1}_{t}$, $\theta_{s} = \theta^{1}_{s} = \theta^{2}_{b}$, $\theta_{b} = \theta^{1}_{b} = \theta^{2}_{s}$, and $\theta_{2} = \theta^{2}_{t}$, which completes the proof of the forward implication. 

The backward implication is a direct consequence of Kron reduction. Starting from a unique solution to the original system, we can calculate the power injections of the reduced systems using~\eqref{eq:kron_reduction}. For area $i\in\mathcal{A}$, let $p^{i}$ be the reduced power injections, $p^{\alpha^{i}}$ be the actual power injection prior to Kron reduction, and $p^{\beta^i}$ be the power injections of the eliminated buses. The power injection at $n\in\mathcal{R}$ such that $n\in\mathcal{R}^{i}_{s}\cap\mathcal{R}^{j}_{b}$, i.e., bus $n$ is inside area~$i$ and shared with area $j$, is
\begin{align}
    &p^{i}_{n} = p^{\alpha^i}_{n} + \sum_{v\in\beta^{i}} A^{i}_{nv}~p^{\beta^i}_{v},\nonumber\\
    &p^{j}_{n} = p^{\alpha^j}_{n} + \sum_{v\in\beta^{j}} A^{j}_{nv}~p^{\beta^j}_{v}.\nonumber
\end{align}
Since both reduced systems have bus $n$, we equate the actual value of the power injection at bus $n$ prior to Kron reduction, i.e., $p^{\alpha^{i}}_{n} = p^{\alpha^{j}}_{n}$. Thus, we have
\begin{subequations}
\begin{align}
& p^{i}_{n} - \sum_{v\in\beta^{i}} A^{i}_{nv}~p^{\beta^{i}}_{v} = p^{j}_{n} -  \sum_{v\in\beta^{j}} A^{j}_{nv}~p^{\beta^{j}}_{v}, \\
& p^{i}_{n} - \sum_{a\in\mathcal{A}} \sum_{v\in\beta^{i}\cap\mathcal{N}^{a}} \!\!\!\!\!\! A^{i}_{nv}~p^{a}_{v} = p^{j}_{n} -  \sum_{a\in\mathcal{A}} \sum_{v\in\beta^{j}\cap\mathcal{N}^{a}} \!\!\!\!\!\! A^{j}_{nv}~p^{a}_{v}, \label{eq:consistency_2} \\
& p^{i}_{n} - \sum_{v\in\beta^{i}\cap\mathcal{N}^{j}} \!\!\!\!\!\! A^{i}_{nv}~p^{j}_{v} - \sum_{a\in\mathcal{A}_{-j}} \sum_{v\in\beta^{i}\cap\mathcal{N}^{a}}  \!\!\!\!\!\! A^{i}_{nv}~p^{a}_{v} = \nonumber\\
& \qquad\qquad p^{j}_{n} - \sum_{v\in\beta^{j}\cap\mathcal{N}^{i}} \!\!\!\!\!\! A^{j}_{nv} p^{i}_{v} - \sum_{a\in\mathcal{A}_{-i}} \sum_{v\in\beta^{j}\cap\mathcal{N}^{a}} \!\!\!\!\!\! A^{j}_{nv}~p^{a}_{v}.\label{eq:consistency_3}
\end{align}
\end{subequations}
\noindent Equation~\eqref{eq:consistency_2} distributes $\beta^{i}$ and $\beta^{j}$ to the buses inside the other areas defined by the sets $\mathcal{N}^{a}$ for $a\in\mathcal{A}$, and~\eqref{eq:consistency_3} pulls the terms that contain variables of the buses inside area $i$ and $j$ out of the summation. Rearranging the terms, we obtain the consistency constraints~\eqref{eq:dist_consistency}, which completes the proof.
\end{proof}

\section{Test Systems and Numerical Results} \label{ap:results}

The list of the test systems, the number of components, and the numerical results of the PTDF and phase-angle formulations are shown in Table~\ref{tab:PTDF_comparison}. The optimal objective column is produced by solving central DCOPF problems. The relative gap is the percentage relative change in the objective function of the distributed method with respect to the central solutions. The time column is the wall-clock time measured in seconds, and the number of iterations is denoted as ``Itr.''. We denote non-convergent instances with ``NC''.
\begin{table*}[b]
\caption{Simulation results of the Phase-Angle and the proposed PTDF Formulations using 44 Test Cases from {PGLib-OPF}}
\begin{tabular}{lrrrrrrrrrrr|}
\cline{7-12}
\multicolumn{1}{l}{} & \multicolumn{1}{l}{} &\multicolumn{1}{l}{} & \multicolumn{1}{l}{} & \multicolumn{1}{l}{} & \multicolumn{1}{l|}{} & \multicolumn{3}{|c|}{\textbf{Phase-Angle DCOPF}}& \multicolumn{3}{|c|}{\textbf{Proposed PTDF-OPF}} \\ \hline
\multicolumn{1}{|l|}{\textbf{Case}} & \multicolumn{1}{|c|}{\textbf{Bus}} & \multicolumn{1}{|c|}{\textbf{Line}} & \multicolumn{1}{|c|}{\textbf{Gen.}} & \multicolumn{1}{|c|}{\textbf{Area}} & \multicolumn{1}{|c|}{\textbf{\begin{tabular}[c]{@{}c@{}}Optimal \\ Objective (\$)\end{tabular}}} & \multicolumn{1}{|c|}{\textbf{\begin{tabular}[c]{@{}c@{}}Relative \\ Gap (\%)\end{tabular}}} & \multicolumn{1}{c|}{\textbf{\begin{tabular}[c]{@{}c@{}}Time \\ (sec.)\end{tabular}}} & \multicolumn{1}{c|}{\textbf{Itr.}} & \multicolumn{1}{|c|}{\textbf{\begin{tabular}[c]{@{}c@{}}Relative \\ Gap (\%)\end{tabular}}} & \multicolumn{1}{c|}{\textbf{\begin{tabular}[c]{@{}c@{}}Time \\ (sec.)\end{tabular}}} & \textbf{Itr.} \\ \hline\hline
\multicolumn{1}{|l|}{\textbf{case57\_ieee}}  & \multicolumn{1}{|r|}{$57$} & \multicolumn{1}{|r|}{$80$} & \multicolumn{1}{|r|}{$7$} & \multicolumn{1}{|r|}{$5$} & \multicolumn{1}{|r|}{$34772$} & \multicolumn{1}{r|}{$9.23\times 10^{-3}$} & \multicolumn{1}{r|}{$2.52$} & \multicolumn{1}{|r|}{$371$} & \multicolumn{1}{r|}{$1.71\times 10^{-2}$} & \multicolumn{1}{r|}{$1.12$} & \multicolumn{1}{r|}{$122$} \\ \hline 
\multicolumn{1}{|l|}{\textbf{case60\_c}} & \multicolumn{1}{|r|}{$60$} & \multicolumn{1}{|r|}{$88$}   & \multicolumn{1}{|r|}{$23$} & \multicolumn{1}{|r|}{$5$} & \multicolumn{1}{|r|}{$90700$} & \multicolumn{1}{r|}{$4.66 \times 10^{-5}$}  & \multicolumn{1}{r|}{14.76} & \multicolumn{1}{r|}{$1875$} & \multicolumn{1}{|r|}{$4.83\times 10^{-4}$} & \multicolumn{1}{r|}{$0.90$}  & \multicolumn{1}{r|}{$121$} \\ \hline
\multicolumn{1}{|l|}{\textbf{case73\_ieee\_rts}} & \multicolumn{1}{|r|}{$73$} & \multicolumn{1}{|r|}{$120$} & \multicolumn{1}{|r|}{$99$} & \multicolumn{1}{|r|}{$3$} & \multicolumn{1}{|r|}{$183003$} & \multicolumn{1}{r|}{$1.38\times 10^{-4}$} & \multicolumn{1}{r|}{$1.88$} & \multicolumn{1}{|r|}{$312$} & \multicolumn{1}{r|}{$1.29\times 10^{-3}$} & \multicolumn{1}{r|}{$0.60$} & \multicolumn{1}{r|}{$65$} \\ \hline
\multicolumn{1}{|l|}{\textbf{case89\_pegase}} & \multicolumn{1}{|r|}{$89$} & \multicolumn{1}{|r|}{$210$} & \multicolumn{1}{|r|}{$12$} & \multicolumn{1}{|r|}{$5$} & \multicolumn{1}{|r|}{$105044$} & \multicolumn{1}{r|}{$1.82\times 10^{-1}$} & \multicolumn{1}{r|}{$336.46$} & \multicolumn{1}{|r|}{$33940$} & \multicolumn{1}{r|}{$1.15\times 10^{-1}$} & \multicolumn{1}{r|}{$3.11$} & \multicolumn{1}{r|}{$371$} \\ \hline
\multicolumn{1}{|l|}{\textbf{case118\_ieee}} & \multicolumn{1}{|r|}{$118$} & \multicolumn{1}{|r|}{$186$} & \multicolumn{1}{|r|}{$54$} & \multicolumn{1}{|r|}{$5$} & \multicolumn{1}{|r|}{$93101$} & \multicolumn{1}{r|}{$1.93\times 10^{-3}$} & \multicolumn{1}{r|}{$5.42$} & \multicolumn{1}{|r|}{$655$} & \multicolumn{1}{r|}{$2.77\times 10^{-4}$} & \multicolumn{1}{r|}{$1.65$} & \multicolumn{1}{r|}{$179$} \\ \hline
\multicolumn{1}{|l|}{\textbf{case162\_ieee\_dtc}} & \multicolumn{1}{|r|}{$162$} & \multicolumn{1}{|r|}{$284$} & \multicolumn{1}{|r|}{$12$} & \multicolumn{1}{|r|}{$5$} & \multicolumn{1}{|r|}{$101462$} & \multicolumn{1}{r|}{$2
65\times 10^{-3}$} & \multicolumn{1}{r|}{$22.08$} & \multicolumn{1}{|r|}{$2250$} & \multicolumn{1}{r|}{$1.03\times 10^{-3}$} & \multicolumn{1}{r|}{$1.96$} & \multicolumn{1}{r|}{$226$} \\ \hline
\multicolumn{1}{|l|}{\textbf{case179\_goc}} & \multicolumn{1}{|r|}{$179$} & \multicolumn{1}{|r|}{$263$} & \multicolumn{1}{|r|}{$29$} & \multicolumn{1}{|r|}{$3$} & \multicolumn{1}{|r|}{$751881$} & \multicolumn{1}{r|}{$2.99\times 10^{-2}$} & \multicolumn{1}{r|}{$0.86$} & \multicolumn{1}{|r|}{$110$} & \multicolumn{1}{r|}{$1.81\times 10^{-1}$} & \multicolumn{1}{r|}{$0.50$} & \multicolumn{1}{r|}{$69$} \\ \hline
\multicolumn{1}{|l|}{\textbf{case197\_snem}} & \multicolumn{1}{|r|}{$197$} & \multicolumn{1}{|r|}{$286$} & \multicolumn{1}{|r|}{$35$} & \multicolumn{1}{|r|}{$5$} & \multicolumn{1}{|r|}{$1.47$} & \multicolumn{1}{r|}{$2.04\times 10^{-3}$} & \multicolumn{1}{r|}{1.77} & \multicolumn{1}{|r|}{200} & \multicolumn{1}{r|}{$1.52\times 10^{-4}$} & \multicolumn{1}{r|}{$0.62$} & \multicolumn{1}{r|}{63} \\ \hline
\multicolumn{1}{|l|}{\textbf{case200\_activ}} & \multicolumn{1}{|r|}{$200$} & \multicolumn{1}{|r|}{$245$} & \multicolumn{1}{|r|}{$49$} & \multicolumn{1}{|r|}{$5$} & \multicolumn{1}{|r|}{$27480$} & \multicolumn{1}{r|}{$1.14\times 10^{-3}$} & \multicolumn{1}{r|}{$10.38$} & \multicolumn{1}{|r|}{$1191$} & \multicolumn{1}{r|}{$4.70\times 10^{-5}$} & \multicolumn{1}{r|}{$0.75$} & \multicolumn{1}{r|}{$87$} \\ \hline
\multicolumn{1}{|l|}{\textbf{case300\_ieee}} & \multicolumn{1}{|r|}{$300$} & \multicolumn{1}{|r|}{$411$} & \multicolumn{1}{|r|}{$69$} & \multicolumn{1}{|r|}{$5$} & \multicolumn{1}{|r|}{$517851$} & \multicolumn{1}{r|}{$7.48\times 10^{-3}$} & \multicolumn{1}{r|}{$37.52$} & \multicolumn{1}{|r|}{$3253$} & \multicolumn{1}{r|}{$7.96\times 10^{-3}$} & \multicolumn{1}{r|}{$2.31$} & \multicolumn{1}{r|}{$220$} \\ \hline
\multicolumn{1}{|l|}{\textbf{case500\_goc}} & \multicolumn{1}{|r|}{$500$} & \multicolumn{1}{|r|}{$733$} & \multicolumn{1}{|r|}{$224$} & \multicolumn{1}{|r|}{$5$} & \multicolumn{1}{|r|}{$440549$} & \multicolumn{1}{r|}{$9.73\times 10^{-5}$} & \multicolumn{1}{r|}{$54.02$} & \multicolumn{1}{|r|}{$4126$} & \multicolumn{1}{r|}{$9.52\times 10^{-5}$} & \multicolumn{1}{r|}{$2.10$} & \multicolumn{1}{r|}{$91$} \\ \hline
\multicolumn{1}{|l|}{\textbf{case588\_sdet}} & \multicolumn{1}{|r|}{$588$} & \multicolumn{1}{|r|}{$686$} & \multicolumn{1}{|r|}{$167$} & \multicolumn{1}{|r|}{$8$} & \multicolumn{1}{|r|}{$310126$} & \multicolumn{1}{r|}{$3.97\times 10^{-3}$} & \multicolumn{1}{r|}{$94.33$} & \multicolumn{1}{|r|}{$7375$} & \multicolumn{1}{r|}{$8.13\times 10^{-4}$} & \multicolumn{1}{r|}{$15.13$} & \multicolumn{1}{r|}{$978$} \\ \hline
\multicolumn{1}{|l|}{\textbf{case793\_goc}} & \multicolumn{1}{|r|}{$793$} & \multicolumn{1}{|r|}{$913$} & \multicolumn{1}{|r|}{$214$} & \multicolumn{1}{|r|}{$5$} & \multicolumn{1}{|r|}{$258308$} & \multicolumn{1}{r|}{$2.08\times 10^{-3}$} & \multicolumn{1}{r|}{$23.79$} & \multicolumn{1}{|r|}{$1603$} & \multicolumn{1}{r|}{$8.06\times 10^{-4}$} & \multicolumn{1}{r|}{$5.44$} & \multicolumn{1}{r|}{$330$} \\ \hline
\multicolumn{1}{|l|}{\textbf{case1354\_pegase}} & \multicolumn{1}{|r|}{$1354$} & \multicolumn{1}{|r|}{$1991$} & \multicolumn{1}{|r|}{$260$} & \multicolumn{1}{|r|}{$5$} & \multicolumn{1}{|r|}{$1218182$} & \multicolumn{1}{r|}{$3.33\times 10^{-3}$} & \multicolumn{1}{r|}{$726.24$} & \multicolumn{1}{|r|}{$35906$} & \multicolumn{1}{r|}{$1.05\times 10^{-3}$} & \multicolumn{1}{r|}{$30.48$} & \multicolumn{1}{r|}{$499$} \\ \hline 
\multicolumn{1}{|l|}{\textbf{case1803\_snem}} & \multicolumn{1}{|r|}{$1803$} & \multicolumn{1}{|r|}{$2795$} & \multicolumn{1}{|r|}{$230$} & \multicolumn{1}{|r|}{4} & \multicolumn{1}{|r|}{$87696$} & \multicolumn{1}{r|}{$3.43\times 10^{-1}$} & \multicolumn{1}{r|}{$4.85$} & \multicolumn{1}{|r|}{$162$} & \multicolumn{1}{r|}{$3.17\times 10^{-1}$} & \multicolumn{1}{r|}{$10.75$} & \multicolumn{1}{r|}{$259$} \\ \hline 
\multicolumn{1}{|l|}{\textbf{case1888\_rte}} & \multicolumn{1}{|r|}{$1888$} & \multicolumn{1}{|r|}{$2531$} & \multicolumn{1}{|r|}{$297$} & \multicolumn{1}{|r|}{$5$} & \multicolumn{1}{|r|}{$1352872$} & \multicolumn{1}{r|}{$2.76\times 10^{-6}$} & \multicolumn{1}{r|}{$1090.74$} & \multicolumn{1}{|r|}{$38265$} & \multicolumn{1}{r|}{$3.10\times 10^{-4}$} & \multicolumn{1}{r|}{$28.59$} & \multicolumn{1}{r|}{$288$} \\ \hline 
\multicolumn{1}{|l|}{\textbf{case1951\_rte}} & \multicolumn{1}{|r|}{$1951$} & \multicolumn{1}{|r|}{$2596$} & \multicolumn{1}{|r|}{$391$} & \multicolumn{1}{|r|}{$5$} & \multicolumn{1}{|r|}{$2031628$} & \multicolumn{1}{r|}{$2.66\times 10^{-6}$} & \multicolumn{1}{r|}{1000.11} & \multicolumn{1}{|r|}{34724} & \multicolumn{1}{r|}{$1.28\times 10^{-2}$} & \multicolumn{1}{r|}{$21.16$} & \multicolumn{1}{r|}{$196$} \\ \hline 
\multicolumn{1}{|l|}{\textbf{case2000\_goc}} & \multicolumn{1}{|r|}{$2000$} & \multicolumn{1}{|r|}{$3639$} & \multicolumn{1}{|r|}{$384$} & \multicolumn{1}{|r|}{$3$} & \multicolumn{1}{|r|}{$943042$} & \multicolumn{1}{r|}{$1.77\times 10^{-4}$} & \multicolumn{1}{r|}{325.72} & \multicolumn{1}{|r|}{$7191$} & \multicolumn{1}{r|}{$2.08\times 10^{-4}$} & \multicolumn{1}{r|}{$37.08$} & \multicolumn{1}{r|}{$264$} \\ \hline 
\multicolumn{1}{|l|}{\textbf{case2312\_goc}} & \multicolumn{1}{|r|}{$2312$} & \multicolumn{1}{|r|}{$3013$} & \multicolumn{1}{|r|}{$444$} & \multicolumn{1}{|r|}{$5$} & \multicolumn{1}{|r|}{$440328$} & \multicolumn{1}{r|}{$1.04\times 10^{-4}$} & \multicolumn{1}{r|}{$279.04$} & \multicolumn{1}{|r|}{$7220$} & \multicolumn{1}{r|}{$3.17\times 10^{-2}$} & \multicolumn{1}{r|}{$119.38$} & \multicolumn{1}{r|}{$1320$} \\ \hline 
\multicolumn{1}{|l|}{\textbf{case2383wp\_k}} & \multicolumn{1}{|r|}{$2383$} & \multicolumn{1}{|r|}{$2896$} & \multicolumn{1}{|r|}{$327$} & \multicolumn{1}{|r|}{$4$} & \multicolumn{1}{|r|}{$1804090$} & \multicolumn{1}{r|}{$4.08\times 10^{-4}$} & \multicolumn{1}{r|}{$8.39$} & \multicolumn{1}{|r|}{$49$} & \multicolumn{1}{r|}{$1.36\times 10^{-3}$} & \multicolumn{1}{r|}{$109.92$} & \multicolumn{1}{r|}{$183$} \\ \hline 
\multicolumn{1}{|l|}{\textbf{case2736sp\_k}} & \multicolumn{1}{|r|}{$2736$} & \multicolumn{1}{|r|}{$3504$} & \multicolumn{1}{|r|}{$420$} & \multicolumn{1}{|r|}{$4$} & \multicolumn{1}{|r|}{$1276034$} & \multicolumn{1}{r|}{$6.09 \times 10^{-6}$} & \multicolumn{1}{r|}{$15.37$} & \multicolumn{1}{|r|}{$113$} & \multicolumn{1}{r|}{$1.46\times 10^{-3}$} & \multicolumn{1}{r|}{$94.49$} & \multicolumn{1}{r|}{$178$} \\ \hline 
\multicolumn{1}{|l|}{\textbf{case2737sop\_k}} & \multicolumn{1}{|r|}{$2737$} & \multicolumn{1}{|r|}{$3506$} & \multicolumn{1}{|r|}{$399$} & \multicolumn{1}{|r|}{$4$} & \multicolumn{1}{|r|}{$764009$} & \multicolumn{1}{r|}{$1.44\times 10^{-4}$} & \multicolumn{1}{r|}{$12.60$} & \multicolumn{1}{|r|}{$104$} & \multicolumn{1}{r|}{$1.84\times 10^{-3}$} & \multicolumn{1}{r|}{$64.66$} & \multicolumn{1}{r|}{$161$} \\ \hline
\multicolumn{1}{|l|}{\textbf{case2742\_goc}} & \multicolumn{1}{|r|}{$2742$} & \multicolumn{1}{|r|}{$4673$} & \multicolumn{1}{|r|}{$182$} & \multicolumn{1}{|r|}{$5$} & \multicolumn{1}{|r|}{$259696$} & \multicolumn{1}{r|}{$1.62\times 10^{-4}$} & \multicolumn{1}{r|}{$296.80$} & \multicolumn{1}{|r|}{$7921$} & \multicolumn{1}{r|}{$3.33\times 10^{-3}$} & \multicolumn{1}{r|}{$20.10$} & \multicolumn{1}{r|}{$169$} \\\hline 
\multicolumn{1}{|l|}{\textbf{case2746wop\_k}} & \multicolumn{1}{|r|}{$2746$} & \multicolumn{1}{|r|}{$3514$} & \multicolumn{1}{|r|}{$520$} & \multicolumn{1}{|r|}{$4$} & \multicolumn{1}{|r|}{$1581425$} & \multicolumn{1}{r|}{$8.68\times 10^{-4}$} & \multicolumn{1}{r|}{$9.63$} & \multicolumn{1}{|r|}{$55$} & \multicolumn{1}{r|}{$8.62\times 10^{-4}$} & \multicolumn{1}{r|}{$183.87$} & \multicolumn{1}{r|}{$139$} \\\hline 
\multicolumn{1}{|l|}{\textbf{case2848\_rte}} & \multicolumn{1}{|r|}{$2848$} & \multicolumn{1}{|r|}{$3776$} & \multicolumn{1}{|r|}{$547$} & \multicolumn{1}{|r|}{$5$} & \multicolumn{1}{|r|}{$1267732$} & \multicolumn{1}{r|}{$8.04\times 10^{-6}$} & \multicolumn{1}{r|}{$1229.30$} & \multicolumn{1}{|r|}{$33304$} & \multicolumn{1}{r|}{$1.80\times 10^{-2}$} & \multicolumn{1}{r|}{$36.79$} & \multicolumn{1}{r|}{$166$} \\\hline 
\multicolumn{1}{|l|}{\textbf{case2853\_sdet}} & \multicolumn{1}{|r|}{$2853$} & \multicolumn{1}{|r|}{$3921$} & \multicolumn{1}{|r|}{$946$} & \multicolumn{1}{|r|}{$5$} & \multicolumn{1}{|r|}{$2036958$} & \multicolumn{1}{r|}{$1.85\times 10^{-4}$} & \multicolumn{1}{r|}{$737.56$} & \multicolumn{1}{|r|}{$16012$} & \multicolumn{1}{r|}{$3.44\times 10^{-4}$} & \multicolumn{1}{r|}{$509.39$} & \multicolumn{1}{r|}{$768$} \\\hline 
\multicolumn{1}{|l|}{\textbf{case2868\_rte}} & \multicolumn{1}{|r|}{$2868$} & \multicolumn{1}{|r|}{$3808$} & \multicolumn{1}{|r|}{$599$} & \multicolumn{1}{|r|}{$5$} & \multicolumn{1}{|r|}{$1966684$} & \multicolumn{1}{r|}{$4.86\times 10^{-6}$} & \multicolumn{1}{r|}{$895.64$} & \multicolumn{1}{|r|}{$25805$} & \multicolumn{1}{r|}{$8.87\times 10^{-2}$} & \multicolumn{1}{r|}{$50.68$} & \multicolumn{1}{r|}{$197$} \\ \hline
\multicolumn{1}{|l|}{\textbf{case2869\_pegase}} & \multicolumn{1}{|r|}{2869} & \multicolumn{1}{|r|}{4582} & \multicolumn{1}{|r|}{510} & \multicolumn{1}{|r|}{$5$} & \multicolumn{1}{|r|}{2386379} & \multicolumn{1}{r|}{\text{NC}} & \multicolumn{1}{r|}{\text{NC}} & \multicolumn{1}{r|}{\text{NC}} & \multicolumn{1}{r|}{$1.51\times 10^{-1}$} & \multicolumn{1}{r|}{$67.82$} & \multicolumn{1}{r|}{$317$} \\ \hline
\multicolumn{1}{|l|}{\textbf{case3012wp\_k}} & \multicolumn{1}{|r|}{$3012$} & \multicolumn{1}{|r|}{$3572$} & \multicolumn{1}{|r|}{$502$} & \multicolumn{1}{|r|}{$2$} & \multicolumn{1}{|r|}{$2509001$} & \multicolumn{1}{r|}{$1.00\times 10^{-5}$} & \multicolumn{1}{r|}{$38.23$} & \multicolumn{1}{|r|}{$263$} & \multicolumn{1}{r|}{$5.00\times 10^{-4}$} & \multicolumn{1}{r|}{$131.31$} & \multicolumn{1}{r|}{$96$} \\\hline
\multicolumn{1}{|l|}{\textbf{case3022\_goc}} & \multicolumn{1}{|r|}{$3022$} & \multicolumn{1}{|r|}{$4135$} & \multicolumn{1}{|r|}{$637$} & \multicolumn{1}{|r|}{$5$} & \multicolumn{1}{|r|}{$599221$} & \multicolumn{1}{r|}{$3.78\times 10^{-4}$} & \multicolumn{1}{r|}{$184.83$} & \multicolumn{1}{|r|}{$4262$} & \multicolumn{1}{r|}{$9.61\times 10^{-2}$} & \multicolumn{1}{r|}{$68.08$} & \multicolumn{1}{r|}{$467$} \\\hline
\multicolumn{1}{|l|}{\textbf{case3120sp\_k}} & \multicolumn{1}{|r|}{$3120$} & \multicolumn{1}{|r|}{$3693$} & \multicolumn{1}{|r|}{$505$} & \multicolumn{1}{|r|}{$2$} & \multicolumn{1}{|r|}{$2087975$} & \multicolumn{1}{r|}{$1.20\times 10^{-3}$} & \multicolumn{1}{r|}{$41.03$} & \multicolumn{1}{|r|}{$291$} & \multicolumn{1}{r|}{$9.63\times 10^{-5}$} & \multicolumn{1}{r|}{$98.18$} & \multicolumn{1}{r|}{$111$} \\\hline
\multicolumn{1}{|l|}{\textbf{case3375wp\_k}} & \multicolumn{1}{|r|}{$3374$} & \multicolumn{1}{|r|}{$4161$} & \multicolumn{1}{|r|}{$596$} & \multicolumn{1}{|r|}{$2$} & \multicolumn{1}{|r|}{$7317012$} & \multicolumn{1}{r|}{$5.98\times 10^{-4}$} & \multicolumn{1}{r|}{$2243.34$} & \multicolumn{1}{|r|}{$18025$} & \multicolumn{1}{r|}{$1.09\times 10^{-4}$} & \multicolumn{1}{r|}{$144.10$} & \multicolumn{1}{r|}{$78$} \\ \hline 
\multicolumn{1}{|l|}{\textbf{case3970\_goc}} & \multicolumn{1}{|r|}{$3970$} & \multicolumn{1}{|r|}{$6641$} & \multicolumn{1}{|r|}{$383$} & \multicolumn{1}{|r|}{$5$} & \multicolumn{1}{|r|}{$934219$} & \multicolumn{1}{r|}{$9.88\times 10^{-5}$} & \multicolumn{1}{r|}{$490.49$} & \multicolumn{1}{|r|}{$9397$} & \multicolumn{1}{r|}{$1.62\times 10^{-4}$} & \multicolumn{1}{r|}{$62.50$} & \multicolumn{1}{r|}{$131$} \\ \hline 
\multicolumn{1}{|l|}{\textbf{case4020\_goc}} & \multicolumn{1}{|r|}{$4020$} & \multicolumn{1}{|r|}{$6988$} & \multicolumn{1}{|r|}{$352$} & \multicolumn{1}{|r|}{$5$} & \multicolumn{1}{|r|}{$795062$} & \multicolumn{1}{r|}{\text{NC}} & \multicolumn{1}{r|}{\text{NC}} & \multicolumn{1}{|r|}{\text{NC}} & \multicolumn{1}{r|}{$4.67\times 10^{-4}$} & \multicolumn{1}{r|}{$144.10$} & \multicolumn{1}{r|}{$485$} \\ \hline 
\multicolumn{1}{|l|}{\textbf{case4601\_goc}} & \multicolumn{1}{|r|}{$4601$} & \multicolumn{1}{|r|}{$7199$} & \multicolumn{1}{|r|}{$408$} & \multicolumn{1}{|r|}{$5$} & \multicolumn{1}{|r|}{$793814$} & \multicolumn{1}{r|}{$2.65\times 10^{-4}$} & \multicolumn{1}{r|}{$466.21$} & \multicolumn{1}{|r|}{$9543$} & \multicolumn{1}{r|}{$4.99\times 10^{-5}$} & \multicolumn{1}{r|}{$62.99$} & \multicolumn{1}{r|}{$118$} \\ \hline
\multicolumn{1}{|l|}{\textbf{case4619\_goc}} & \multicolumn{1}{|r|}{$4619$} & \multicolumn{1}{|r|}{$8150$} & \multicolumn{1}{|r|}{$347$} & \multicolumn{1}{|r|}{$5$} & \multicolumn{1}{|r|}{$457436$} & \multicolumn{1}{r|}{$5.53\times 10^{-5}$} & \multicolumn{1}{r|}{$2056.33$} & \multicolumn{1}{|r|}{$32988$} & \multicolumn{1}{r|}{$4.13\times 10^{-3}$} & \multicolumn{1}{r|}{$54.57$} & \multicolumn{1}{r|}{$207$} \\ \hline
\multicolumn{1}{|l|}{\textbf{case4661\_sdet}} & \multicolumn{1}{|r|}{$4661$} & \multicolumn{1}{|r|}{$5997$} & \multicolumn{1}{|r|}{$1176$} & \multicolumn{1}{|r|}{$5$} & \multicolumn{1}{|r|}{$2216303$} & \multicolumn{1}{r|}{$5.69\times 10^{-4}$} & \multicolumn{1}{r|}{$1017.50$} & \multicolumn{1}{|r|}{$12210$} & \multicolumn{1}{r|}{$7.87\times 10^{-5}$} & \multicolumn{1}{r|}{$738.49$} & \multicolumn{1}{r|}{$841$} \\ \hline
\multicolumn{1}{|l|}{\textbf{case4837\_goc}} & \multicolumn{1}{|r|}{$4837$} & \multicolumn{1}{|r|}{$7765$} & \multicolumn{1}{|r|}{$332$} & \multicolumn{1}{|r|}{$5$} & \multicolumn{1}{|r|}{$850397$} & \multicolumn{1}{r|}{$2.05\times 10^{-4}$} & \multicolumn{1}{r|}{$1166.83$} & \multicolumn{1}{|r|}{$18926$} & \multicolumn{1}{r|}{$2.35\times 10^{-3}$} & \multicolumn{1}{r|}{$49.16$} & \multicolumn{1}{r|}{$235$} \\ \hline
\multicolumn{1}{|l|}{\textbf{case4917\_goc}} & \multicolumn{1}{|r|}{$4917$} & \multicolumn{1}{|r|}{$6726$} & \multicolumn{1}{|r|}{$1349$} & \multicolumn{1}{|r|}{$5$} & \multicolumn{1}{|r|}{$1383655$} & \multicolumn{1}{r|}{$2.20\times 10^{-3}$} & \multicolumn{1}{r|}{$826.88$} & \multicolumn{1}{|r|}{$10759$} & \multicolumn{1}{r|}{$5.07\times 10^{-2}$} & \multicolumn{1}{r|}{$366.19$} & \multicolumn{1}{r|}{$857$} \\ \hline
\multicolumn{1}{|l|}{\textbf{case5658\_epigrids}} & \multicolumn{1}{|r|}{$5658$} & \multicolumn{1}{|r|}{$9078$} & \multicolumn{1}{|r|}{$474$} & \multicolumn{1}{|r|}{$5$} & \multicolumn{1}{|r|}{$1195466$} & \multicolumn{1}{r|}{$4.47\times 10^{-5}$} & \multicolumn{1}{r|}{$3310.68$} & \multicolumn{1}{|r|}{$52147$} & \multicolumn{1}{r|}{$2.25\times 10^{-3}$} & \multicolumn{1}{r|}{$103.51$} & \multicolumn{1}{r|}{$217$} \\ \hline
\multicolumn{1}{|l|}{\textbf{case6468\_rte}} & \multicolumn{1}{|r|}{$6468$} & \multicolumn{1}{|r|}{$9000$} & \multicolumn{1}{|r|}{$1295$} & \multicolumn{1}{|r|}{$5$} & \multicolumn{1}{|r|}{$1982818$} & \multicolumn{1}{r|}{$2.36\times 10^{-5}$} & \multicolumn{1}{r|}{$1931.97$} & \multicolumn{1}{|r|}{$29198$} & \multicolumn{1}{r|}{$2.20\times 10^{-5}$} & \multicolumn{1}{r|}{$106.24$} & \multicolumn{1}{r|}{$263$} \\ \hline
\multicolumn{1}{|l|}{\textbf{case6470\_rte}} & \multicolumn{1}{|r|}{$6470$} & \multicolumn{1}{|r|}{$9005$} & \multicolumn{1}{|r|}{$1330$} & \multicolumn{1}{|r|}{$5$} & \multicolumn{1}{|r|}{$2136095$} & \multicolumn{1}{r|}{$1.24\times 10^{-4}$} & \multicolumn{1}{r|}{$1870.29$} & \multicolumn{1}{|r|}{$25400$} & \multicolumn{1}{r|}{$4.25\times 10^{-4}$} & \multicolumn{1}{r|}{$335.21$} & \multicolumn{1}{r|}{$517$} \\ \hline
\multicolumn{1}{|l|}{\textbf{case6495\_rte}} & \multicolumn{1}{|r|}{$6495$} & \multicolumn{1}{|r|}{$9019$} & \multicolumn{1}{|r|}{$1372$} & \multicolumn{1}{|r|}{$5$} & \multicolumn{1}{|r|}{$2561787$} & \multicolumn{1}{r|}{$4.00\times 10^{-3}$} & \multicolumn{1}{r|}{$1721.17$} & \multicolumn{1}{|r|}{$23346$} & \multicolumn{1}{r|}{$1.84\times 10^{-4}$} & \multicolumn{1}{r|}{$188.33$} & \multicolumn{1}{r|}{$272$} \\ \hline
\multicolumn{1}{|l|}{\textbf{case6515\_rte}} & \multicolumn{1}{|r|}{$6515$} & \multicolumn{1}{|r|}{$9037$} & \multicolumn{1}{|r|}{$1388$} & \multicolumn{1}{|r|}{$5$} & \multicolumn{1}{|r|}{$2559329$} & \multicolumn{1}{r|}{$9.70\times 10^{-4}$} & \multicolumn{1}{r|}{$2774.59$} & \multicolumn{1}{|r|}{$36943$} & \multicolumn{1}{r|}{$3.30\times 10^{-4}$} & \multicolumn{1}{r|}{$219.62$} & \multicolumn{1}{r|}{$312$} \\
\hline
\end{tabular}\label{tab:PTDF_comparison}
\end{table*}

\balance
\bibliographystyle{IEEEtran}
\bibliography{IEEEabrv,references.bib}

\end{document}